\documentclass[12pt,DIV=11,titlepage=false]{scrartcl}
\usepackage[a4paper,left=2cm,right=1.2cm,top=3cm,bottom=4cm]
{geometry}
\usepackage{amsmath,amstext, amsthm, amsbsy,amssymb,marvosym,fancyhdr,graphicx,amscd,amsfonts,latexsym,delarray,stackrel,lineno,color,cite,latexsym,euscript}
\usepackage{wasysym}
\usepackage{srcltx}
\usepackage{mathrsfs}

\usepackage[all]{xy}
\usepackage{t1enc}

\usepackage{pifont}

\usepackage{mathpple}
\usepackage[T1]{fontenc}

\newcommand*\patchAmsMathEnvironmentForLineno[1]{%
  \expandafter\let\csname old#1\expandafter\endcsname\csname #1\endcsname
  \expandafter\let\csname oldend#1\expandafter\endcsname\csname end#1\endcsname
  \renewenvironment{#1}%
     {\linenomath\csname old#1\endcsname}%
     {\csname oldend#1\endcsname\endlinenomath}}%
\newcommand*\patchBothAmsMathEnvironmentsForLineno[1]{%
  \patchAmsMathEnvironmentForLineno{#1}%
  \patchAmsMathEnvironmentForLineno{#1*}}%
\AtBeginDocument{%
\patchBothAmsMathEnvironmentsForLineno{equation}%
\patchBothAmsMathEnvironmentsForLineno{align}%
\patchBothAmsMathEnvironmentsForLineno{flalign}%
\patchBothAmsMathEnvironmentsForLineno{alignat}%
\patchBothAmsMathEnvironmentsForLineno{gather}%
\patchBothAmsMathEnvironmentsForLineno{multline}%
}

\definecolor{Green}{rgb}{0,1,0}
\definecolor{Blue}{RGB}{0,0,191}
\definecolor{mathmodecolor}{RGB}{0,102,0}
\definecolor{keywordcolor}{RGB}{0,51,151}
\definecolor{sourcebackgroundcolor}{RGB}{255,247,223}
\definecolor{unixagred}{RGB}{255,0,0}
\definecolor{lightgray}{RGB}{191,191,191}
\definecolor{green}{RGB}{1,191,191}

\newtheorem{thm}{Theorem}[section]
\newtheorem{prop}[thm]{Proposition}
\newtheorem{cor}[thm]{Corollary}
\newtheorem{lem}[thm]{Lemma}

\newtheorem{defn}[thm]{Definition}
\newtheorem{rem}[thm]{Remark}
\newtheorem{example}[thm]{Example}

\def\End{{\rm End}}

\def\Hom{{\rm Hom}}
\def\id{{\rm id}}

\def\sign{{\rm sign}}

\def\Spec{{\rm Spec\,}}
\def\Sp{{\rm Spec}\,}

\def\sign{\mathbf{S}}

\def\B{{\mathbb B}}
\def\C{{\mathbb C}}
\def\F{{\mathbb F}}
\def\K{{\mathbb K}}
\def\N{{\mathbb N}}

\def\Q{{\mathbb Q}}
\def\R{{\mathbb R}}
\def\Z{{\mathbb Z}}

\def\supp{{\rm Support}}

\def\hatz{{\hat\Z^\times}}

\def\vvert{{\Vert}}

\def\spzb{{\overline{\Spec\Z}}}

\def\coprod{{\amalg}}

\def\cA{{\mathcal A}}
\def\cB{{\mathcal B}}
\def\cC{{\mathcal C}}

\def\cG{{\mathcal G}}
\def\cH{{\mathcal H}}

\def\cM{{\mathcal M}}
\def\cO{{\mathcal O}}

\def\cR{{\mathcal R}}
\def\cS{{\mathcal S}}

\def\gam{{\Gamma\Ses}}
\def\ssei{{\underline\Hom_\gop}}

\def\qqq{\,,\,~\forall}

\newcommand{\ie}{{\it i.e.\/}\ }
\newcommand{\eg}{{\it e.g.\/}\ }
\newcommand{\cf}{{\it cf.}}
\newcommand{\opcit}{{\it op.cit.\/}\ }
\newcommand{\resp}{{\it resp.,\/}\ }

\def\hq{{\cO_\spzb}}
\def\hqr{{{\cO}_\epsilon}}

\def\spz{{\Spec\Z}}

\def\id{{\mbox{Id}}}

\def\Hom {{\mbox{Hom}}}

\def\End{{\mbox{End}}}

\def\fii{{\rm  finite}}
\def\ssr{{\mathfrak{  SR}}}

\def\sss{{\mathbb  S}}
\def\gop{{\Gamma^{\rm op}}}

\def\sign{\mathbf{S}}

\def\Se{\frak{ Sets}}
\def\fin{\frak{ Fin}}

\newcommand{\nil}[1]{}

\def\Ses{{\Se_*}}
\def\sses{{\cS_*}}

\def\mono{{\rm Monoid}}
\def\salg{{\sss{\rm alg}}}
\def\spzb{{\overline{\Spec\Z}}}
\begin{document}

    \title{Absolute algebra and Segal's $\Gamma$-rings}
    \subtitle{\LARGE au dessous de $\overline{\Sp(\Z)}$}
    \author{~\\Alain~\textsc{Connes\,}\thanks{Coll\`ege de France,
3 rue d'Ulm, Paris F-75005 France. 
I.H.E.S. and Ohio State University. ~email: alain@connes.org}~ ~and ~Caterina~\textsc{Consani\,}\thanks{Department of Mathematics,  Johns Hopkins
University, Baltimore MD 21218 USA.~ email: kc@math.jhu.edu \newline The authors are grateful to the referees for very helpful comments.}} 

\date{}

    \maketitle

\vspace{-.4in}

\begin{abstract}
\noindent{\bf \large Abstract} ~~We show that the basic categorical concept of an $\sss$-algebra as derived from the theory of Segal's $\Gamma$-sets provides a unifying description of several constructions  attempting to model an algebraic geometry over the absolute point. It merges, in particular, the approaches using mono\"{\i}ds, semirings and hyperrings as well as the development by means of monads and generalized rings in Arakelov geometry.  The assembly map  determines a functorial way to associate an $\sss$-algebra to a monad on pointed sets.  The 
notion of an $\sss$-algebra is very familiar in algebraic topology where it also provides a suitable groundwork to the definition of topological cyclic homology. The main contribution of  this paper is to point out its relevance and unifying role in arithmetic,  in relation with the development of an algebraic geometry over symmetric closed monoidal categories. 
\end{abstract}

\tableofcontents

\section{Introduction}

 The notion of an $\sss$-algebra (\ie algebra over the sphere spectrum)  is well-known in homotopy theory (\cf~\eg\cite{DGM}):  in the categorical form used in this paper and implemented by the concept of a discrete $\Gamma$-ring, it was formalized in the late 90's. $\sss$-algebras are also intimately related to the theory of brave new rings (first introduced in the 80's) and of functors with smash products (FSP) for spectra in algebraic topology (\cf~\cite{Bok,Lyd0,Lyd}). 
The goal of this article is to explain how the implementation of the notion of an $\sss$-algebra in arithmetic, in terms of Segal's $\Gamma$-rings, succeeds to unify several constructions pursued in recent times aimed to define the notion of  ``absolute algebra". In particular, we refer to the development, for applications in number theory and algebraic geometry, of a suitable framework apt to provide
a rigorous meaning  to the process of taking the limit of geometry over finite fields $\F_q$ as $q\to 1$. In our previous work we have met and used at least three possible categories suitable to handle this unification: namely the category $\mathbf\cM$ of mono\"{\i}ds as in  \cite{CC1,CC3}, the category $\mathbf\cH$ of hyperrings of \cite{CC2,CC4,CC5} and finally the category $\mathbf\cS$ of  semirings as in \cite{C,CCas}.\newline
In \cite{Durov}, N. Durov developed a geometry over $\F_1$ intended for Arakelov theory applications by implementing  monads as  generalizations of classical rings. In his work certain combinatorial structures replace, at the archimedean place(s), the geometric constructions classically performed, at a finite prime ideal,  in the process of reduction modulo that ideal. In the present  article we argue that all the relevant constructions in \opcit can be naturally subsumed by the well-known theory of $\sss$-algebras and Segal's $\Gamma$-sets\footnote{The more general notion of $\Gamma$-space is needed in homotopy theory and is simply that of a simplicial $\Gamma$-set.} in homotopy theory. This theory  is taken up as the natural groundwork in the recent book \cite{DGM}. In Proposition \ref{assembly} we prove  that the assembly map of \cite{Lyd} provides a functorial way to associate an $\sss$-algebra  to a monad on pointed sets.  While in the context of  \cite{Durov} the tensor product $\Z\otimes_{\F_1} \Z$ produces an uninteresting output isomorphic to $\Z$, in Proposition \ref{tensq} we show that the same tensor square, re-understood in the theory of $\sss$-algebras, provides a highly non-trivial object. As explained in \cite{Lyd}, the category  of pointed $\Gamma$-sets  is a symmetric closed monoidal category and the theory of generalized schemes of B. T\"oen and M. Vaqui\'e in \cite{TV}  applies directly to this category while it could not be implemented in the category of endofunctors under composition of \cite{Durov} which is not symmetric. \newline
We endow the Arakelov compactification $\spzb$ of $\Spec\Z$ with a natural structure of a  sheaf $\cO_\spzb$ of $\sss$-algebras and each Arakelov divisor provides a natural sheaf of modules over that structure sheaf. Moreover, this new structure of $\spzb$ over $\sss$ endorses a one parameter group of weakly invertible sheaves whose tensor product rule is the same as   the composition rule of the Frobenius correspondences over
the arithmetic site \cite{CCas, CCas1}.    
The fundamental advantage of  having unified the various attempts done in recent times to provide a suitable definition of  ``absolute algebra" by means of  the well established concept of $\sss$-algebra is  that this latter notion is at the root of  the theory of topological cyclic homology which can be understood as cyclic homology over the absolute base $\sss$. Segal's $\Gamma$-spaces  form a model for stable homotopy theory and  they may be equivalently viewed as simplicial objects in the category of $\Gamma$-sets (\cf\cite{DGM}). Note that general $\Gamma$-spaces in \opcit are not assumed to be special or very special, and it is only by dropping the very special condition that one can use the Day product.  Moreover,
the prolongation of a $\Gamma$-space to a spectrum is a purely formal
construction and does not need the special/very special hypotheses.  The very special condition only appears for ``fibrant" objects in the appropriate Quillen model category. Thus, the stable homotopy category, understood using general $\Gamma$-spaces,  is properly viewed as the derived category of the  algebraic category of $\Gamma$-sets, \ie of $\sss$-modules. Topological cyclic homology then appears as the direct analogue of ordinary cyclic homology of rings when one replaces the category of abelian groups by the category of $\Gamma$-sets. In particular it is now available to understand the new structure of $\spzb$ using its structure sheaf $\cO_\spzb$ and the modules. \newline
Our original motivation for using cyclic homology in the arithmetic context arose from the following two results:\newline
$(i)$ In \cite{CC6} we showed that cyclic homology (a
fundamental tool in noncommutative geometry)  determines the correct infinite dimensional (co)homological theory for arithmetic varieties to recast the archimedean local factors of Serre as regularized determinants. The key operator in this context is  the generator of the $\lambda$-operations in cyclic theory.\newline
$(ii)$ L. Hesselholt and I. Madsen have proven  that the de Rham-Witt complex, an essential ingredient of crystalline cohomology,  arises 
naturally when one studies the topological cyclic homology of smooth algebras over a perfect field of finite characteristic (\cf~\eg \cite{HM,H}). 

Our long term objective is to use cyclic homology over the absolute base $\sss$ and the arithmetic site defined in \cite{CCas, CCas1} to obtain a global interpretation of $L$-functions of arithmetic varieties.

\section{$\sss$-algebras and Segal's $\Gamma$-sets}\label{sectgamma}

Before recalling the definition of a $\Gamma$-set  we comment on its conceptual yet simple meaning. A $\Gamma$-set is the most embracing generalization of the datum provided on a set by a commutative addition with a zero element. The commutativity of the addition is encoded by the fact  that  the output of a (finite) sum depends only upon the corresponding finite set (with multiplicity). A preliminary definition of a $\Gamma$-set  is given by a covariant  functor $F$ from the category $\fin$ of finite sets   to the category $\Se$ of sets, mapping a finite set $X$ to the set $F(X)$ of all possible strings of elements indexed by $X$ which can be added. The action of $F$ on morphisms encodes the addition, while its covariance  embodies the associativity of the operation. One introduces a base point $\ast\in X$ in the labelling of the sums and a base point $*\in F(X)$ to implement the presence of the $0$-element (neutral for the sum). Thus $F$ is best described  by a pointed covariant functor from the category  $\fin_*$ of finite pointed sets to the category $\Se_*$ of pointed sets. Note that while the category of covariant functors  $\fin_*\longrightarrow\Se$ is a topos the category of pointed covariant functors $\fin_*\longrightarrow\Se_*$
is no longer such and admits an initial object which is also final, in direct analogy with the category of abelian groups. The theory reviewed in Chapter II of \cite{DGM} is that of Segal's $\Gamma$-spaces, \ie of simplicial $\Gamma$-sets. The category of $\Gamma$-spaces admits a natural structure of model category (\cf \opcit Definition 2.2.1.5) allowing one to do homotopical algebra. The transition from $\Gamma$-sets to $\Gamma$-spaces parallels  the transition from abelian groups to chain complexes of abelian groups in positive degrees and  the constructions of topological Hochschild and cyclic homology then become parallel to the construction of their algebraic ancestors.

\subsection{Segal's $\Gamma$-sets}
For each integer $k\geq 0$, one introduces the pointed finite set $k_+:=\{0,\ldots ,k\}$, 
where $0$ is the base point. Let $\Gamma^{\rm op}$ be  the small, full subcategory  of $\fin_*$ whose objects are the sets $k_+$'s, for $k\geq 0$. The morphisms in the category are maps between pointed finite sets, preserving the base points.   The notion  of a discrete $\Gamma$-space, \ie of a $\Gamma$-set  suffices for our applications and at the same time it is intimately related to the topos $\hat \Gamma$ of covariant functors $\Gamma^{\rm op}\longrightarrow \Se$. The small category $\gop$ is a pointed category \ie it admits a (unique) initial and final object, namely the object $0_+$ formed by the base point alone. In general, if $\cC$ is a pointed category one defines a $\Gamma$-object of $\cC$ to be a covariant functor $\gop\longrightarrow \cC$ preserving the base point. 

\begin{defn}\label{defngamset} A $\Gamma$-set $F$ is a  functor $F:\gop\longrightarrow\Ses$ between pointed categories. \end{defn}

The morphisms $\Hom_\gop(M,N)$ between two $\Gamma$-sets are natural transformations of functors.
The category $\gam$ of $\Gamma$-sets is a symmetric, closed, monoidal category (\cf~\cite{DGM}, Chapter II).  The monoidal structure is given by the smash product of $\Gamma$-sets which is an example of a Day product: we shall review this product  in \S \ref{salgsect}. The closed structure property is shown in \cite{Lyd} (\cf~also \cite{DGM} Theorem 2.1.2.4:  replace the category $\sses$ with $\Ses$). The internal $\Hom$ structure is given by 
\begin{equation}\label{gamspec1}
\ssei(M,N):=\{n_+\mapsto \Hom_\gop(M,N_{n_+\wedge})\}
\end{equation} 
where $N_{n_+\wedge}(k_+):=N(n_+\wedge k_+)$. For two pointed sets $X,Y$ their smash product is defined by collapsing the set $X\times \{*\}\cup\{*\}\times Y$ to $\{*\}\times \{*\}$ in the product $X\times Y$.\vspace{.05in}

Our next task is to explain how the closed monoidal category $\gam$ encompasses several attempts to model  a notion of an  ``absolute algebra". 

\subsection{A basic construction of $\Gamma$-sets}
 Let $M$ be a commutative mono\"{\i}d denoted additively and with a $0$ element. One defines a functor $HM:\gop\longrightarrow\Ses$ by setting
\begin{equation}\label{functsum}
HM(k_+):=M^k, \quad  HM(f):HM(k_+)\to HM(n_+), \  \  HM(f)(\phi)(y):=\sum_{x\in f^{-1}(y)} \phi(x).
\end{equation}
In the last formula $y\in \{1,\ldots ,n\}$ and the sum over the empty set is the $0$ element.  Moreover the base point of $HM(k_+):=M^k$ is $\phi(x)=0$, $\forall x\in \{1,\ldots ,k\}$. The formula \eqref{functsum} is meaningful in view of the commutativity and associativity of the  mono\"{\i}d $M$. One easily checks that the maps $f: k_+\to n_+$ with $f^{-1}(\{0\})\neq \{0\}$ do not create any problem for the functoriality. The functor $HM$ is in fact a particular case of the covariant functor $\fin\longrightarrow\Se$, $X\mapsto M^X$ defined by the formula $f_\star\phi(y)=\sum_{x\mid f(x)=y} \phi(x)$. More precisely,  $HM$ is obtained by restricting to $\fin_*$ and by subsequently dividing using the equivalence relation $$\phi\sim \phi'\iff \phi(x)=\phi'(x) \qquad   \forall x\neq *.$$
By assuming that $f:X\to Y$ preserves the base point and if $\phi\sim \phi'$, one gets  $$\sum_{x\mid f(x)=y} \phi(x)=\sum_{x\mid f(x)=y} \phi'(x), \qquad \forall y\neq *.$$

\subsection{$\sss$-algebras}\label{salgsect}
Next, we recall the notion of an $\sss$-algebra as given in   \cite{DGM} (Definition 2.1.4.1). This requires to define first the smash product of two $\Gamma$-sets, \ie of two pointed functors $F_j: \gop\longrightarrow\Ses$, $j=1,2$. The definition of the smash product $F_1\wedge F_2$ is dictated by the internal $\Hom$ structure of \eqref{gamspec1} and the adjunction formula
$$
\Hom_\gop(F_1\wedge F_2,G)\simeq \Hom_\gop(F_1,\ssei(F_2,G)).
$$
Thus the smash product $F_1\wedge F_2$ is such that for any $\Gamma$-set $G$, a morphism $F_1\wedge F_2\to G$ is simply described by a natural transformation of bifunctors $\gop\times \gop\longrightarrow \Ses$ \ie by maps of pointed sets, natural in both objects $X,Y$ of  $\gop$ 
$$
F_1(X)\wedge F_2(Y)\to G(X\wedge Y).
$$
The evaluation of the $\Gamma$-set $F_1\wedge F_2$  on an object $Z$ of $\fin_*$ is given by the following colimit 
\begin{equation}\label{colim}
(F_1\wedge F_2)(Z)=\varinjlim_{v:X\wedge Y\to Z}( F_1(X)\wedge F_2(Y)), 
\end{equation}
where for any morphisms $f:X\to X'$ and $g:Y\to Y'$ in $\fin_*$ one uses the morphism 
$$
 F_1(f)\wedge  F_2(g): F_1(X)\wedge  F_2(Y)\to F_1(X')\wedge  F_2(Y')
$$
provided that $v'\circ (f\wedge g)=v$, with $v': X'\wedge Y'\to Z$. Thus, with the exception of the base point, a point of $( F_1\wedge  F_2)(Z)$ is represented by the data $(X,Y,v,x,y)$ given by a pair of objects $X,Y$ of $\fin_*$, a map $v:X\wedge Y\to Z$ and a pair of non-base points $x\in  F_1(X)$, $y\in  F_2(Y)$. Moreover, notice the following implication
$$
v'\circ (f\wedge g)=v \implies (X,Y,v,x,y)\sim (X',Y',v', F_1(f)(x), F_2(g)(y)).
$$
The colimit \eqref{colim} is in general not filtered.
The specialization of Definition 2.1.4.1. of \cite{DGM} to the case of $\Gamma$-sets yields the following
\begin{defn}\label{defnsalg} An $\sss$-algebra $\mathcal A$  is a $\Gamma$-set $\mathcal A: \gop\longrightarrow\Ses$  endowed with an associative multiplication 
$\mu:\cA \wedge \cA\to \cA$ and a unit~ $1: \sss\to \cA$, where  
$$
\sss:\gop\longrightarrow \Ses, \qquad  \sss(k_+):=k_+\quad \forall k\geq 0,\quad \sss(f) := f.
$$
is the canonical inclusion    functor.
\end{defn}
Since $\sss$ is the inclusion functor, the obvious identity map  $X\wedge Y\to X\wedge Y$ defines the product $\mu: \sss\wedge\sss\to \sss$ in
$\sss$ which is clearly associative.  By construction, any $\sss$-algebra is an algebra over $\sss$.  This elementary categorical object can be naturally associated to $\F_1$ \ie to the most basic algebraic structure underlying the absolute (geometric) point.

A morphism  $\rho:\cA\to \cB$ of $\sss$-algebras is a morphism of the underlying $\Gamma$-sets which is compatible with the unit and with the product. This latter condition is equivalent to the commutativity of the following diagram for any objects $X,Y$ of $\gop$
\begin{equation}\label{multss}
\xymatrix{
\cA(X)\wedge \cA(Y) \ar[d]^{\mu_{\cA}}  \ar[rr]^{\rho_X\wedge \rho_Y}&& \ar[d]^{\mu_\cB}  \cB(X)\wedge \cB(Y) \\
\cA(X\wedge Y) \ar[rr]_{\rho_{X\wedge Y}}&& \cB(X\wedge Y)
}
\end{equation}
We let $\salg$ be the category of $\sss$-algebras. 
It follows from \opcit (\cf~\S 2.1.4.1.6) that, given an $\sss$-algebra $\cA$ and an integer $n\geq 1$, one obtains an $\sss$-algebra of $n\times n$ matrices  by endowing the $\Gamma$-set $M_n(\cA):=\Hom_\Ses(n_+,n_+\wedge \cA(X))$ with the natural multiplication of matrices having only one non-zero entry in each column. There is moreover a straightforward notion of module over an $\sss$-algebra (\cf~\opcit Definition 2.1.5.1): an $\sss$-module being just a $\Gamma$-set. 

\section{Basic constructions of $\sss$-algebras}

\subsection{$\sss$-algebras and mono\"{\i}ds}\label{sectlevel1}
We begin this section by reviewing an easy construction of an $\sss$-algebra derived from a functor from the category of (not necessarily commutative) multiplicative mono\"{\i}ds with a unit and a zero element, to the category of $\sss$-algebras.  This construction is described in  Example 2.1.4.3, 2. of \cite{DGM}, where the mono\"{\i}d $M$ is not assumed to have a $0$ element and the obtained $\sss$-algebra is called spherical mono\"{\i}d algebra.
\begin{defn} Let $M$ be a multiplicative mono\"{\i}d with a multiplicative unit and a zero element $0$. We define  the covariant functor 
\[
\sss M: \fin_*\longrightarrow \Se_*,  \qquad \sss M(X)=M\wedge X
\]
with $0\in M$ viewed as the base point and with maps $\id_M\times f$, for $f:X\to Y$.
\end{defn}
\begin{prop}\label{mono2sss} Let $M$ be a multiplicative mono\"{\i}d with a unit and a zero element.  Then  the product in $M$  endows $\sss M$ with a structure of an $\sss$-algebra.
\end{prop}
\proof The product in $M$, viewed as a map $M\wedge M\to M$ determines the following map, natural in both objects $X$, $Y$ of $\fin_*$
$$
\sss M(X)\wedge \sss M(Y)=(M\wedge X)\wedge (M\wedge Y)\to M\wedge X\wedge Y=\sss M(X\wedge Y).
$$
The multiplicative unit $1\in M$ determines a natural transformation $\sss \to \sss M$. This construction endows $\sss M$ with a structure of an $\sss$-algebra.\endproof
The multiplicative mono\"{\i}d $\{0,1\}$ (frequently denoted by $\F_1$) determines, in this framework, the canonical inclusion functor $\sss$: \ie $\sss\{0,1\} = \sss$.

Conversely, given an $\sss$-algebra  $\cA:\gop\longrightarrow \Ses$  one obtains, using the product $\mu: \cA(1_+)\wedge \cA(1_+)\to \cA(1_+)$ and the unit $1_\cA:\sss\longrightarrow \cA$,   a canonical structure of multiplicative mono\"{\i}d  (with a base point  $0$ and a multiplicative identity $1$) on the set $M=\cA(1_+)$.\begin{prop}\label{monomono} Let $\cA$ be an $\sss$-algebra, $M$ a mono\"{\i}d and $\alpha:M\to \cA(1_+)$ a morphism of  mono\"{\i}ds, then the following map defines a morphism of $\sss$-algebras
$$
\tilde \alpha:\sss M\longrightarrow \cA,\quad \tilde \alpha(j\wedge m)=\mu_\cA(1_\cA(j)\wedge \alpha(m))\in \cA(k_+\wedge 1_+)=\cA(k_+)\qqq j\in k_+,~ m\in M.
$$
The map $\alpha\mapsto \tilde \alpha$ gives the adjunction $\Hom_\mono(M,\cA(1_+))\simeq \Hom_\salg(\sss M,\cA)$.
\end{prop}
\proof The unit $1_\cA:\sss\longrightarrow \cA$ for $\cA$ defines a natural transformation of functors compatible with the product. Since the product $\mu: \cA(X)\wedge \cA(Y)\to \cA(X\wedge Y)$ is natural in $X$ and $Y$, this shows that by taking $Y=1_+$, the morphism $\tilde \alpha$ defines a natural transformation $\tilde \alpha:\sss M\longrightarrow \cA$. The associativity of $\mu$ shows that $\tilde \alpha$ is multiplicative. The required adjunction then follows using \eqref{multss}.\endproof
Note that the counit of the adjunction gives a canonical morphism $\rho:\sss\cA(1_+)\to \cA$.

\subsection{From semirings to $\sss$-algebras}\label{hyperpart}

A semiring $R$ is a set endowed with two binary operations: $+$ and $\cdot$. The addition $+$ defines on $R$ the structure of an additive, commutative mono\" id with neutral element $0\in R$; the multiplication $x,y\mapsto xy$ is left and right distributive with respect to the addition and defines on $R$ the structure of a multiplicative
mono\" id with neutral element $1$ and absorbing element $0$.  We let $\ssr$ be the category of semirings.
An important construction of an $\sss$-algebra  is provided by the following result (\cf\cite{DGM})
\begin{lem}\label{sssalg} Let $R$ be a semiring, then the functor $HR: \gop \to \Ses$, $X\mapsto HR(X)= R^{X/*}$ is naturally endowed with a structure of an $\sss$-algebra.
\end{lem}
\proof The additive structure on $R$ defines the $\Gamma$-set $HR$ using \eqref{functsum}. We first describe the product and the unit $1:\sss\to HR$. To define the product we introduce the following map, natural in both objects $X,Y$ of  $\gop$ 
$$
HR(X)\wedge HR(Y)\to HR(X\wedge Y), \  \  (\phi,\psi)\mapsto \phi\psi, \  \  \phi\psi(x,y)=\phi(x)\psi(y)\qqq x\in X\setminus \{*\}, y\in Y\setminus \{*\}.
$$
The naturality of the above operation follows from the bilinearity of the product in $R$. The unit $1:\sss\to HR$ is defined as follows 
$$
1_X:X\to HR(X), \  \  1_X(x)=\delta_x, \ \  \delta_x(y):= \begin{cases} 0 & \text{if} \ x\neq y
\\
1 & \text{if} \ x=y.
\end{cases}
$$
One obtains in this way  a natural transformation since there is at most one non-zero value in the sum $\sum_{x\mid f(x)=y} \delta_a(x)$ which computes $HR(f)(\delta_a)$. Moreover, a non-zero value occurs exactly when $y=f(a)$, which shows that $HR(f)(\delta_a)=\delta_{f(a)}$. Notice that we have defined both the product and the transformation $1:\sss\to HR$ without using the additive group structure of a ring: in fact the semiring structure suffices. Finally, the axioms of $\sss$-algebras are checked in the same way as for rings (\cf \cite{DGM} Example 2.1.4.3 for details). \endproof

\begin{prop}\label{sssalg1} The functor $H: \ssr \to \salg$~ from semirings to $\sss$-algebras is fully faithful.
\end{prop}
\proof One needs to check  that, for two semirings $A,B$ the natural map 
$$
\Hom_{\ssr}(A,B)\to \Hom_{\salg}(HA,HB), \  \ \alpha \mapsto H\alpha
$$
is bijective. As a set, one obtains the bijection $HA(1_+)\stackrel{\epsilon}{\simeq}A$ by using the map defined by $\epsilon(\phi)=\phi(1)\in A$, $\forall \phi \in 
HA(1_+)$, which maps the base point to $0\in A$. The product in $A$ is recovered by the map $HA(1_+)\wedge HA(1_+)\to HA(1_+\wedge 1_+)$, using the fact that $1_+\wedge 1_+=1_+$.\newline
Notice that for a pointed functor $\cA:\gop \longrightarrow \Ses$ of the form $\cA=HA$, with $A$ a semiring, one derives the special property that given two elements $x,y\in HA(1_+)$ there exists a {\em unique} element $z\in HA(2_+)$ whose images by the maps $\alpha,\beta :2_+\to 1_+$ of the form 
\begin{equation}\label{alphabeta}
\alpha:0\mapsto 0, \ 1\mapsto 1, \ 2\mapsto 0, \  \  \  \beta: 0\mapsto 0, \ 1\mapsto 0,\ 2\mapsto 1
\end{equation}
are given by $HA(\alpha)z=x$,  $HA(\beta)z=y$. One then obtains  
\begin{equation}\label{addition}
x+y=\epsilon(HA(\gamma)z),\quad\text{for}\quad \gamma:0\mapsto 0, \ 1\mapsto 1, \ 2\mapsto 1.
\end{equation}
This shows that a morphism of functors $\rho\in\Hom_{\salg}(HA,HB)$ determines, by restriction to $1_+$, a homomorphism $\rho_1\in \Hom_{\ssr}(A,B)$ of semirings. The  uniqueness of this homomorphism is clear using the bijection $\epsilon:HA(1_+)\to A$. The equality $\rho=H(\rho_1)$ follows from the naturality of $\rho$ which implies that  the projections $\epsilon_j^A:HA(k_+)\to HA(1_+)=A$, $\epsilon_j^A(\phi)=\phi(j)$, fulfill
$$
\rho(\phi)(j)=\epsilon_j^B(\rho(\phi))=\rho_1(\epsilon_j^A(\phi))=\rho_1(\phi(j)).
$$
\endproof 

The formula \eqref{addition} for the addition  in $HA(1_+)$ retains  a meaning also in the case where the above special property,  for a pointed functor $\cA:\gop \longrightarrow \Ses$ is relaxed by dropping the condition that the solution to  $\cA(\alpha)z=x$,  $\cA(\beta)z=y$ is unique. In this case, one can still  define the following generalized addition \ie the hyper-operation which associates to a pair $x,y\in \cA(1_+)$ a subset of $ \cA(1_+)$
\begin{equation}\label{hypersum}
x \oplus y:=\{\cA(\gamma)z\mid z\in \cA(2_+), \  \cA(\alpha)z=x, \ \cA(\beta)z=y\}\qqq x,y\in \cA(1_+).
\end{equation}
This fact suggests that one can  associate an $\sss$-algebra to  a hyperring. This construction will be described in more details in  Section  \ref{secthyp}.

\section{Smash products}\label{sectsmashprod}

The apparent simplicity of the smash product of two $\Gamma$-sets in the Day product  hides in fact  a significant difficulty of the concrete computation of the colimit defined in \eqref{colim}.
In this section we consider  the  specific example of the $\sss$-algebra $H\B$, where $\B:=\{0,1\},\,  1+1=1$ is the smallest semiring
of characteristic one (idempotent) and compute explicitly the smash product $H\B\wedge H\B$. We start by  producing an explicit description of the $\sss$-algebra $H\B$. 

\subsection{The $\sss$-algebra $H\B$}

\begin{lem}\label{sssalg1} The $\sss$-algebra $H\B$ is canonically isomorphic to the functor $P:\gop \longrightarrow \Ses$ which associates to an object $X$ of $\gop$ the set $P(X)$ of subsets of $X$ containing the base point. The functoriality is given by the direct image $X\supset A\mapsto f(A)\subset Y$. The smash product is provided by the map $P(X)\wedge P(Y)\to P(X\wedge Y)$ which associates to the pair $(A,B)$ the smash product $A\wedge B\subset X\wedge Y$. The unit  $1:\sss\to H\B$ is given by the natural transformation $X\to P(X)$ defined by the map  $X\ni x\mapsto \{*,x\}\subset X$.
\end{lem}
\proof A map $\phi:X\to \B$ is specified by the subset $\phi^{-1}(\{1\})=\{x\in X\mid \phi(x)=1\}$. One associates to $\phi\in H\B(X)$ the subset $A=\{*\}\cup \phi^{-1}(\{1\})\subset X$. The formula $\sum_{x\mid f(x)=y} \phi(x)$ which defines the functoriality shows that it corresponds to the direct image $X\supset A\mapsto f(A)\subset Y$. The last two statements are straightforward to check using the product in $\B$.  \qed

\begin{rem}\label{f2}{\rm It is interesting to compare the structures of the $\sss$-algebras $H\B$ and $H\F_2$ because, when evaluated on an object $X$ of $\gop$, they both yield the same set. $H\F_2(X)$ is indeed equal to the set $P(X)$ of subsets of $X$ containing the base point: this equality is provided  by associating to  $\phi\in H\F_2(X)$ the subset $A=\{*\}\cup \phi^{-1}(\{1\})\subset X$. The products and the unit maps   are the same in both constructions since  the product in $\F_2$ is the same as that in $\B$. The difference between the two constructions becomes finally visible by analyzing  the functoriality property of the maps. The simple rule $X\supset A\mapsto f(A)\subset Y$ holding in $H\B$ is replaced, in the case of $H\F_2$ by adding the further condition that on the set $f(A)\subset Y$ one only retains the points $y\in f(A)$ such that the cardinality of the preimage $f^{-1}(\{y\})$ is an odd number.}
\end{rem}

\subsection{The set $(H\B\wedge H\B)(k_+)$}

It follows from the colimit definition \eqref{colim} that $(P\wedge P)(k_+)$ is the set $\pi_0(\cC_k)$ of connected components of the following category $\cC_k$.
The  objects of $\cC_k$ are 4-tuples $\alpha=(X,Y,v,E)$, 
where $X,Y$ are objects of $\fin_*$, $v:X\wedge Y\to k_+$ is a morphism in that category and  $E \in P(X)\wedge P(Y)$. A morphism $\alpha\to \alpha'$ in $\cC_k$ is given by a  pair of pointed maps $f:X\to X'$, $g:Y\to Y'$ such that $v'\circ (f,g)=v$ and $P(f)\wedge P(g)(E)=E'$.
 
Notice that the set $P(X)$  of subsets of $X$ containing the base point can be equivalently viewed as the set of all subsets of $X\setminus \{*\}$. Thus   $E\in P(X)\wedge P(Y)$, $E\neq *\in P(X)\wedge P(Y)$, is  determined by a pair  $(A,B)$ of non-empty subsets,
$A\subset X\setminus \{*\}$, $B\subset Y\setminus \{*\}$. One can thus encode an object $\alpha=(X,Y,v,E)$ of $\cC_k$ such that $E\neq *$ by a 5-tuple  $\alpha=(X,Y,v,(A,B))$ where $A\subset X\setminus \{*\}$, $B\subset Y\setminus \{*\}$ are non-empty sets.  A morphism $\alpha\to \alpha'$ of such 5-tuples is then  given by a pair of pointed maps  $f:X\to X', \ g:Y\to Y'$ such that 
\begin{equation}\label{equiv1} f(A_+)=A'_+,\ g(B_+)=B'_+,\ v'\circ (f,g)=v.
\end{equation}
For $\alpha'=(X',Y',v',*)$, a morphism $\alpha\to \alpha'$ is given by a pair of pointed maps  $f:X\to X', \ g:Y\to Y'$ such that
\begin{equation}\label{equiv1bis}
 v'\circ (f,g)=v, \  \  f(A_+)=\{*\},\ \text{or}\ g(B_+)=\{*\}
\end{equation}
\begin{defn}
	Let $\alpha=(X,Y,v,E)$ be an object of $\cC_k$ then the {\em support} of $\alpha$ is empty if $E=*$ and for $E\neq *\in P(X)\wedge P(Y)$,  $E=(A,B)$ one lets
	\begin{equation}\label{equiv1} 
	\supp(\alpha):=(A\times B)\cap W, \ \ W:=v^{-1}(k_+\setminus \{0\})\subset (X\setminus \{*\})\times (Y\setminus \{*\})
	\end{equation}
	An object $\alpha$ of $\cC_k$ is {\em degenerate} iff $\supp(\alpha)=\emptyset$.
\end{defn}
In particular one sees that any object $\alpha=(X,Y,v,*)$ is degenerate.
\begin{lem}\label{assoc0}
Let $(f,g):\alpha\to \alpha'$ be a morphism in $\cC_k$ then 
\begin{equation}\label{equiv12} 
	(f\times g)(\supp(\alpha))=\supp(\alpha')	\end{equation}
\end{lem}
\proof Assume first that $\alpha'=(X',Y',v',*)$. Then by \eqref{equiv1bis} the restriction of $v$ to $A_+\wedge B_+$ is equal to $0$. Thus $\supp(\alpha)=\emptyset$. Assume now that $\alpha'=(X',Y',v',(A',B'))$.
For $(x,y)\in (A\times B)\cap W$, one has $v'(f(x),g(y))=v(x,y)\neq 0$, hence $f(x)\in X'\setminus \{*\}$, $g(y)\in Y'\setminus \{*\}$ and $(f(x),g(y))\in (A'\times B')\cap W'$. Conversely,  let $(x',y')\in (A'\times B')\cap W'$, then since $f(A_+)=A'_+,\ g(B_+)=B'_+$ there exist $x\in  A$, $f(x)=x'$, $y\in  B$, $g(y)=y'$.
One then has $(x,y)\in(A\times B)\cap W$ since $v(x,y)=v'(f(x),g(y))=v'(x',y')\neq 0$. This shows that $(f\times g)(\supp(\alpha))=\supp(\alpha')$.\endproof
Let $\star:=(*,*,\id,*)$ be the object of $\cC_k$ where $X$ and $Y$ are reduced to the base point. One has 
\begin{lem}\label{assoc05}
An object of $\cC_k$ is degenerate if and only if it belongs to the connected  component of $\star$. 
\end{lem}
\proof The condition $\supp(\alpha)=\emptyset$ defines a connected family of objects by Lemma \ref{assoc0}. Thus it is enough to show that if $\supp(\alpha)=\emptyset$ then $\alpha$ is in the same component as $\star$. 
The inclusions 
$
\iota_A:A_+\subset X, \ \iota_B:B_+\subset Y 
$
give in general  a morphism in $\cC_k$
\begin{equation}\label{firstred}
(\iota_A,\iota_B):(A_+,B_+,v\vert_{(A_+\wedge B_+)},A,B)\to (X,Y,v,A,B)	
\end{equation}
If $\supp(\alpha)=\emptyset$ the restriction of $v$ to $A_+\wedge B_+$ is equal to $0$. This shows that the unique pair of maps $p:A_+\to *$, $q:B_+\to *$ gives a morphism in $\cC_k$
$$
(p,q):(A_+,B_+,v\vert_{(A_+\wedge B_+)},A,B)\to (*,*,\id,*).
$$
so that $\alpha=(X,Y,v,A,B)$ and $\star$ are in the same component.
\endproof 

\subsection{k-relations}
To determine an explicit description of the (pointed) set $(H\B\wedge H\B)(k_+)$ we introduce the following terminology
\begin{defn}\label{kmultigraph} $(i)$~A $k$-{\em relation}   is a triple $C=(F,G,v)$ where $F$ and $G$ are non-empty finite sets and $v:F\times G\to k_+$ is a map of sets such that no line or column of the corresponding matrix is identically $0$\footnote{\ie for each $x\in F$ there exists $y\in G$ such that $v(x,y)\neq 0$ and symmetrically $\forall y\in G$, $\exists x\in F\mid v(x,y)\neq 0$.}.  

$(ii)$~A morphism $C\to C'$ between two k-relations $C=(F,G,v)$ and $C'=(F',G',v')$ 
is a pair of surjective maps $f:F\to F'$, $g:G\to G'$ such that $v'\circ (f\times g)=v$.

$(iii)$~A $k$-relation~$C=(F,G,v)$ is said to be reduced when no line is repeated and no column is repeated.\end{defn}
We let $\cG_k$ be the category of $k$-relations. The next lemma produces a retraction of the
full subcategory of $\cC_k$ whose objects are non-degenerate,  on $\cG_k$. 
\begin{lem}\label{assoc}$(i)$~Let $\alpha=(X,Y,v,A,B)$ be a non-degenerate  object of $\cC_k$. Set $F=p_A(\supp(\alpha))$,  $G=p_B(\supp(\alpha))$
where $p_A$ and $p_B$ are the projections. Then the triple 
$
\Gamma(\alpha)=(F,G,v\vert_{F\times G})
$
defines a k-relation.

$(ii)$~A morphism $\rho:\alpha\to \alpha'$ in $\cC_k$ induces by restriction a morphism $\Gamma(\rho):\Gamma(\alpha)\to \Gamma(\alpha')$ of k-relations.

$(iii)$~For $\gamma=(F,G,v)$ a k-relation, let  $\tilde \gamma=(F_+,G_+,\tilde v,F,G)$, where $\tilde v:F_+\wedge G_+\to k_+$ is the unique extension of $v$ as a pointed map. Then $\gamma\mapsto \tilde \gamma$ extends to a functor from the category $\cG_k$ of $k$-relations to $\cC_k$ and one has $\Gamma(\tilde \gamma)=\gamma$.

$(iv)$~The objects  $\alpha$ and $ \tilde \Gamma(\alpha)$ belong to the same connected component of $\cC_k$.
\end{lem}
\proof $(i)$~Notice that $(A\times B)\cap W$ is a subset of $A\times B$ whose projections are $(p_A((A\times B)\cap W)=F$ and $p_B((A\times B)\cap W)=G$, thus  the restriction of $w$ to $F\times G$ defines a k-relation. \newline
$(ii)$~Let $(f,g):\alpha\to \alpha'$ be a morphism in $\cC_k$ between non-degenerate objects. One has a commutative diagram  
\begin{equation}\label{comdi1}
\xymatrix{
\supp(\alpha) \ar[d]^{p_A}  \ar[rr]^{f\times g}&& \ar[d]^{p_{A'}}  \supp(\alpha') \\
F \ar[rr]_{f}&& F'
}
\end{equation}
where both $f\times g$ and $p_{A'}$ are surjective.
 This shows that $f(p_A(\supp(\alpha)))=p_{A'}(\supp(\alpha'))$ and similarly that $g(p_B(\supp(\alpha)))=p_{B'}(\supp(\alpha'))$. Thus the restrictions of $f$ and $g$ define surjections $f:F\to F'$ and $g:G\to G'$, where $F=p_A(\supp(\alpha))$, $G=p_B(\supp(\alpha))$. The restrictions $v\vert_{F\times G}:F\times G\to k_+$, $v'\vert_{F'\times G'}:F'\times G'\to k_+$ fulfill the equality $v=v'\circ (f,g)$ and one thus obtains  a morphism $\Gamma(\alpha)\to \Gamma(\alpha')$ of k-relations.\newline
$(iii)$~Let $(f,g):\gamma=(F,G,v)\to (F',G',v')=\gamma'$ be a morphism of k-relations, we extend both $f,g$ to maps of pointed sets $f_+:F_+\to F'_+$, $g_+:G_+\to G'_+$. Then one easily checks that one obtains a morphism $\tilde \gamma\to \tilde \gamma'$ (\cf \eqref{equiv1}). One has by construction $\Gamma(\tilde \gamma)=\gamma$.\newline 
$(iv)$~By \eqref{firstred} the object $\alpha=(X,Y,v,A,B)$ is equivalent to $\beta=(A_+,B_+,v\vert_{(A_+\wedge B_+)},A,B)$. 
We now prove that $\beta$ is equivalent to $ \tilde \Gamma(\alpha)=(F_+,G_+,\tilde w,F,G)$, where  $F=p_A((A\times B)\cap W)$, $G=p_B((A\times B)\cap W)$, and $w=v\vert_{(F\times G)}$. 
Let  $f:A_+\to F_+$ act as the identity on $F\subset A$ and map the elements of $A\setminus F$ to the base point. We define $g:B_+\to G_+$ in a similar manner. Then by construction one has $f(A_+)=F_+,\ g(B_+)=G_+$. Moreover for $(a,b)\in A\times B$ one has $w(a,b)=v(f(a),g(b))$ since both sides vanish unless $(a,b)\in F\times G$ and they agree on $F\times G$. This shows that $\alpha$ and $ \tilde \Gamma(\alpha)$ belong to the same connected component. \endproof

Every k-relation $C=(F,G,v)$ admits a canonical reduction $r(C)=(F/\!\!\sim,~ G/\!\!\sim,~ v\vert_{(F/\sim\times G/\sim)})$. This is obtained by dividing $F$ and $G$ by the following equivalence relations  
$$
i\sim i'\iff v(i,j)=v(i',j)\ \forall j\in G, \ \ \ j\sim j'\iff v(i,j)=v(i,j')\ \forall i\in F.
$$
By construction, the value $v(i,j)$ only depends upon the classes $(x,y)$ of $(i,j)$  and thus the canonical reduction map $r_C:C\to r(C)$ defines a morphism of k-relations. More precisely
\begin{lem}\label{assoc17}$(i)$~The reduction map defines an endofunctor $r:\cG_k\longrightarrow \cG_k$ of the category of $k$-relations which determines a retraction on the full subcategory of reduced $k$-relations.\newline
$(ii)$~Let $(f,g): C\to C'$ be a morphism of $k$-relations, with $C$ reduced,  then $f$ and $g$ are bijections and $(f,g)$ is an isomorphism.
\end{lem}
\proof $(i)$~Let $(f,g):C=(F,G,v)\to C'=(F',G',v')$ be a morphism of $k$-relations.
We show that there exists a unique morphism $(r(f),r(g)):r(C)\to r(C')$ making the following diagram commutative
\begin{equation}\label{redrel127}
\xymatrix{
C \ar[d]^{r_C}  \ar[rr]^{(f,g)}&& \ar[d]^{r_{C'}}  C' \\
r(C) \ar[rr]_{(r(f),r(g))}&& r(C')
}
\end{equation} 
Indeed, for $i,i'\in F$ one has 
$
i\sim i'\implies f(i)\sim f(i')
$, since using the equality $v'\circ (f,g)=v$ and the surjectivity of $g$ one derives
$$
v(i,j)=v(i',j)\ \forall j\in G \implies v'(f(i),j')=v'(f(i'),j')\ \forall j'\in G'.
$$
Similarly for $j,j'\in G$  one has 
$
j\sim j'\implies g(j)\sim g(j')
$. Thus $f$ and $g$ induce  maps $r(f)$ and $r(g)$ on the reductions making the 
diagram \eqref{redrel127} commutative. The surjectivity of $f$ and $g$ implies  the surjectivity of the induced maps $r(f)$ and $r(g)$.

$(ii)$~Let $(f,g):C=(F,G,v)\to C'=(F',G',v')$ be a morphism of $k$-relations, with $C$ reduced. Then for $i\neq i'\in F$ let $j\in G$ such that $v(i,j)\neq v(i',j)$. One has 
$$
v'(f(i),g(j))=v(i,j)\neq v(i',j)=v'(f(i'),g(j))
$$
which shows that $f$ is injective. Thus both $f$ and $g$ are bijections and $(f,g)$ is an isomorphism. \endproof

Let $\cR(k)$ be the set of isomorphism classes of reduced $k$-relations, and let $\cR(k)_+$ be the set obtained by adjoining a base point. We extend the association $k_+\mapsto \cR(k)_+$ to a functor  $\cR_+:\gop\longrightarrow\Se_*$ by letting a morphism $\phi\in\Hom_\gop(k_+,\ell_+)$ act on elements of $\cR(k)_+$ as follows 
$$
\cR_+(\phi)(F,G,v):=r\circ \Gamma((F,G,\phi\circ v))
$$
where the right hand side reduces to the base point when $(F,G,\phi\circ v)$ is degenerate.
\begin{thm}\label{wedge2} $(i)$~For  $k\in \N$, the map $r\circ \Gamma$ induces a canonical bijection of sets
 $$
 r\circ \Gamma:(H\B\wedge H\B)(k_+)\to \cR(k)_+
 $$

$(ii)$~The isomorphism $ r\circ \Gamma$ induces an equivalence of functors $H\B\wedge H\B\simeq \cR_+$.
\end{thm}
\proof $(i)$~By construction, $(P\wedge P)(k_+)$ is the union of its base point with  the set $\pi_0(\cC_k)$ of connected components of the category $\cC_k$. It follows from Lemma \ref{assoc} that the category $\cG_k$ of $k$-relations is a retraction of the category of non-degenerate objects of $\cC_k$. Thus using Lemma \ref{assoc05} one obtains $\pi_0(\cC_k)\simeq \pi_0(\cG_k)_+$.
Moreover by Lemma \ref{assoc17} the reduction $r$ gives a bijection 
$\pi_0(\cG_k)\simeq \cR(k)$.
\newline
$(ii)$~This follows since by construction of $H\B\wedge H\B$ one has
$$
H\B\wedge H\B(\phi)(X,Y,v,E)=(X,Y,\phi\circ v,E)
$$
which can be applied to any representative $(X,Y,v,E)$ of a given class in $\pi_0(\cC_k)$.
\endproof 

\begin{cor}\label{pp1} $(i)$~For $n\in \N$, let $\id_n$ be the graph of the identity map on the set with $n$-elements. Then the $1$-relations $\id_n$ belong to distinct connected components of the category $\cG_1$ and they define distinct elements of $(H\B\wedge H\B)(1_+)$.

$(ii)$~The action of the cyclic group $C_2$ on $(H\B\wedge H\B)(1_+)$ is the transposition acting on isomorphism classes of reduced $1$-relations. Both its set of fixed points and its complement are infinite.\end{cor}
\proof $(i)$~The statement follows from the fact that the $1$-relations $\id_n$ are reduced and pairwise non-isomorphic.\newline
$(ii)$~The action of the transposition $\sigma \in C_2$ replaces a reduced 1-relation $\alpha$ by its transpose $\alpha^t$. The fixed points are given by the reduced 1-relations such that $\alpha$ is isomorphic to $\alpha^t$;  all the $\id_n$ as in $(i)$ are thus fixed points. The relations between finite sets of different cardinality determine infinitely many non fixed points. \endproof

\begin{figure}
\begin{center}
\includegraphics[scale=0.8]{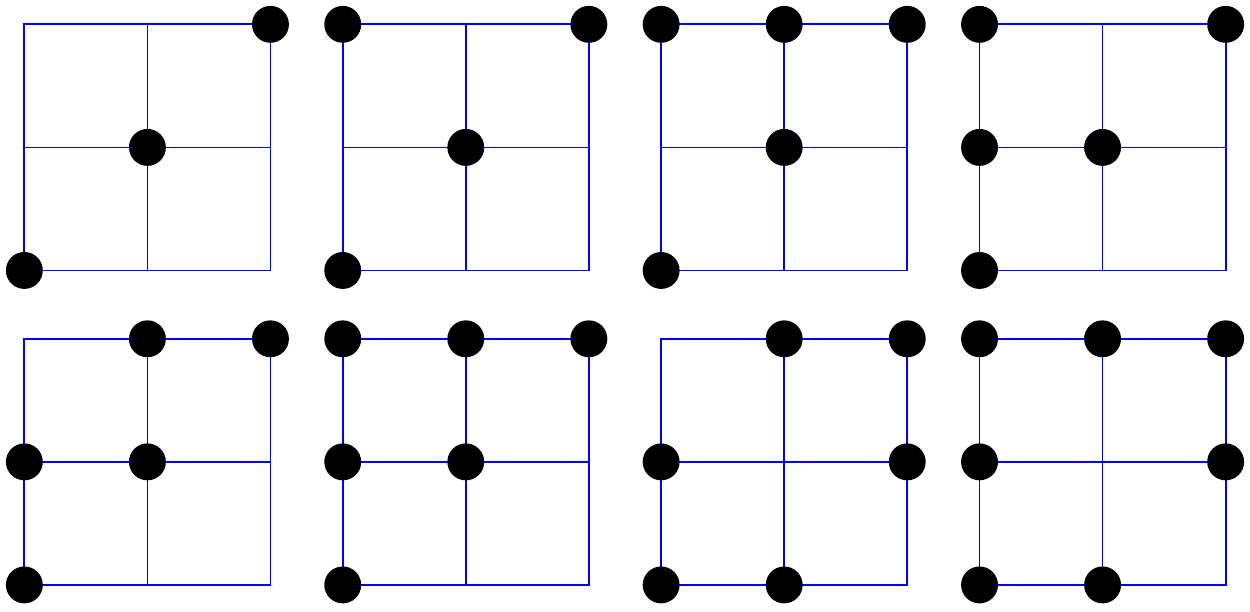}
\end{center}
\caption{Reduced $1$-relations $(3,3)$\label{krelations9} }
\end{figure}

\begin{rem}{\rm $(i)$~In general, the natural map $M(1_+)\wedge N(1_+)\to (M\wedge N)(1_+)$ is {\em not} surjective. This can be easily seen by choosing $M=N=H\B$, since $H\B(1_+)$ is finite unlike $(H\B\wedge H\B)(1_+)$.

$(ii)$~A square reduced relation  is not necessarily symmetric. For instance, among the $8$ reduced $1$-relations  of Figure \ref{krelations9}, only the following reduced 1-relation and its transposed
$$ \left( \begin{array}{ccc}
1 & 1 & 1 \\
1 & 0 & 0 \\
0 & 1 &0 \end{array} \right), \qquad \left( \begin{array}{ccc}
1 & 1 & 0 \\
1 & 0 & 1 \\
1 & 0 &0 \end{array} \right)
$$
 which  have equal number (3) of lines and columns are  {\em not symmetric},  
since the number of (non-zero) elements in the lines are respectively $(3,1,1)$ and $(2,2,1)$ which cannot be matched by any pair of permutations. 
}\end{rem}

\begin{figure}
\begin{center}
\includegraphics[scale=0.5]{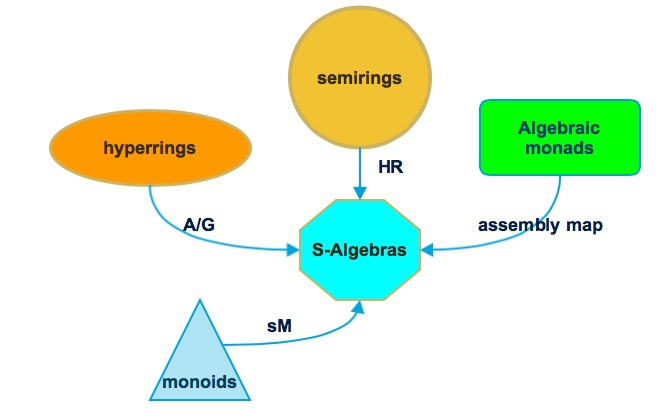}
\end{center}
\caption{Sources of $\sss$-algebras\label{levels} }
\end{figure}

\section{Hyperrings and $\sss$-algebras}\label{secthyp}

The definition \eqref{hypersum} of \S\ref{hyperpart} which gives back  the classical addition on $R$ for the $\sss$-algebra  $HR$ associated to a semiring $R$ is multivalued for a general $\sss$-algebra $\cA: \gop\longrightarrow\Se_*$. The extension of ring theory when the additive group structure is replaced by a multivalued addition has been investigated by M. Krasner in his theory of hyperrings and hyperfields. In our recent research we have encountered these types of generalized algebraic structures in the following two important cases (\cf \cite{CC4,CC5}):

$-$~The ad\`ele class space of a global field $F$ is naturally a hyperring and it contains the Krasner hyperfield $\K=\{0,1\}$ as a sub-hyperalgebra precisely because one divides the ring of the ad\`eles of $F$ by the non-zero elements $F^\times$ of the field $F$.

$-$~The dequantization process described by the  ``universal perfection'' at the archimedean place of a number field yields a natural hyperfield structure $\R^\flat$ on the set of real numbers and  a natural hyperfield structure $\C^\flat$  in the complex case. 

The structure of a hyperfield is more rigid than that of a semifield since the operation $x\mapsto -x$ is present in the hyper-context while it is missing for semirings. This fact makes the requirement that non-zero elements are invertible  more restrictive for hyperfields   than for semifields. A natural construction of  hyperrings is achieved by simply dividing an ordinary (commutative) ring $A$ with a subgroup $G\subset A^\times$ of the group of invertible elements of $A$.  In fact, in Corollary 3.10 of \cite{CC4} we proved that any (commutative) hyperring extension of the Krasner hyperfield $\K=\{0,1\}$ ($1+1 = \{0,1\}$), without zero divisors and of dimension $> 3$ is of the form $A/G$, where $A$ is a commutative ring and $G\subset A$ is the group of non-zero elements of a subfield $F\subset A$. Moreover, Theorem 3.13 of \opcit shows that a morphism of hyperrings of the above form lifts (under mild non-degeneracy conditions) to a morphism of the pairs $(A,F)$.
\subsection{From hyperrings $A/G$ to $\sss$-algebras}\label{fromrtohyp}
In this section we provide the construction of the $\sss$-algebra associated to a hyperring of the form $A/G$, where $A$ is a commutative ring and $G\subset A^\times$ is a subgroup. Notice that the structure needed to define {\em uniquely} the $\sss$-algebra $HA/G$ is provided by the pair $(A,G)$ and that this datum is more precise than simply assigning the hyperring $A/G$.
\begin{prop}\label{sssalg3} Let $A$ be a commutative ring and $G\subset A^\times$ be a subgroup of the group of invertible elements of $A$.
 For each object $X$ of $\fin_*$, let $(HA/G)(X)$
be the quotient of $HA(X)$ by the following equivalence relation 
$$
\phi\sim \psi\iff \exists g\in G~\text{s.t.}~ \psi(x)=g\phi(x)\qqq x\in X, x\neq *.
$$
Then, the functor $HA/G:\gop\longrightarrow\Ses$ defines an $\sss$-algebra and the quotient map $HA\longrightarrow HA/G$ is a morphism of $\sss$-algebras.
\end{prop}
\proof We first check the functoriality of the construction, \ie that for $f:X\to Y$ in $\fin_*$, the map $HA(f)$ respects the equivalence relation $\phi\sim \psi$. For $y\in Y$, $y\neq *$ one has
$$
HA(f)(\psi)(y):=\sum_{x\mid f(x)=y} \psi(x)=\sum_{x\mid f(x)=y} g\phi(x)=g\left(\sum_{x\mid f(x)=y} \phi(x)\right)=gHA(f)(\phi)(y).
$$
Since these equations hold for the same $g$ and for all $y\in Y$, one derives that $HA(f)(\phi)\sim HA(f)(\psi)$. The above equivalence relation is  compatible with the product: in fact using the commutativity of $A$, one has 
$$
(g\phi)(h\psi)(x,y)=g\phi(x)h\psi(y)=gh\phi(x)\psi(y)\qqq x\in X\setminus \{*\}, y\in Y\setminus \{*\}.
$$
The unit $1:\sss\longrightarrow HA$ defines by composition a natural transformation $1:\sss\longrightarrow HA/G$: 
$$
1_X:X\to (HA/G)(X), \ \ 1_X(x)=G\delta_x\in (HA/G)(X)=HA(X)/G.
$$
It follows easily from the construction that the quotient map $HA\longrightarrow HA/G$ is a morphism of $\sss$-algebras.\qed

\begin{prop}\label{sssalgback} Let $A$ be a commutative ring, let $G\subset A^\times$ be a subgroup of the group of invertible elements of $A$ and $HA/G$ the $\sss$-algebra defined in Proposition \ref{sssalg3}. Then, the set $(HA/G)(1_+)$ endowed with the hyper-addition defined in  \eqref{hypersum} and the product $HA/G(1_+)\wedge HA/G(1_+)\to HA/G(1_+)$ induced by the $\sss$-algebra structure,
 is canonically isomorphic to the hyperring $A/G$.
\end{prop}
\proof  By using the bijection $\epsilon(\phi):=\phi(1)\in A/G$, $\forall \phi \in 
(HA/G)(1_+)$ which maps the base point to $0\in A$ one derives the bijection of sets $(HA/G)(1_+)\stackrel{\epsilon}{\simeq}A/G$. The product in $(HA/G)(1_+)$ is recovered by the map $(HA/G)(1_+)\wedge (HA/G)(1_+)\to (HA/G)(1_+\wedge 1_+)$ using the relation $1_+\wedge 1_+=1_+$. It is the same as the product  in $A/G$.
 By using \eqref{alphabeta} and \eqref{addition}, the hyper-addition  \eqref{hypersum} is described, for $x,y\in (HA/G)(1_+)\simeq A/G$ as follows
$$
x \oplus y:=\{(HA/G)(\gamma)z\mid z\in (HA/G)(2_+), \  (HA/G)(\alpha)z=x, \ (HA/G)(\beta)z=y\}
$$
For $x=aG$ and $y=bG$ let $\phi\in HA(2_+)$ with $
\phi(1)=a$ and $\phi(2)=b$. Then the coset $(a+b)G$  belongs to $x \oplus y$ and all elements of $x \oplus y$ are of this form. In this way one recovers the hyper-addition on $A/G$.  \endproof

Given an $\sss$-algebra $\cA$ a necessary condition for  $\cA$ being of the above type $HA/G$ is that the equations $\cA(\alpha)z=x$,  $\cA(\beta)z=y$  always admit solutions. The lack of uniqueness of solutions is directly related to the multivalued nature of the addition.

\subsection{$W$-models and associated $\sss$-algebras}

The construction of the $\sss$-algebras $HA/G$ given in Proposition \ref{sssalg3} applies in particular when one uses the notion of $W$-model introduced in \cite{CC5}. We recall that
given a hyperfield $K$, a $W$-{\em model} of $K$ is by definition a triple $(W,\rho,\tau)$ where\vspace{.05in}

$W$ is a field\newline
$\rho:W\to K$ is a homomorphism of hyperfields\newline
$\tau: K\to W$ is a multiplicative section of $\rho$.\vspace{.05in}

Because $\tau: K\to W$ is a multiplicative section of $\rho$, the map $\rho$ is surjective and identifies the multiplicative mono\"{\i}d $K$  with the quotient of $W$ by the multiplicative subgroup $G=\{x\in W\mid \rho(x)=1\}$. Since $\rho:W\to K$ is a homomorphism of hyperfields, one has $\rho(x+y)\in \rho(x)+\rho(y)$ for any $x,y\in W$. Thus $\rho$ defines a morphism of hyperfields $W/G\to K$.\newline
A morphism $(W,\rho,\tau)\to (W',\rho',\tau')$ of $W$-models is a field homomorphism $\alpha:W\to W'$ such that $\rho'\circ \alpha=\rho$ and $\alpha\circ \tau=\tau'$. We recall that a $W$-model of $K$ is said to be {\em universal}  if it is an initial object in the category of $W$-models of $K$. When such universal model exists one easily sees that it is unique up to canonical isomorphism: we denote it by $W(K)$. In this case it is thus natural to associate to $K$ the $\sss$-algebra $H(W(K))/G$ of its universal $W$-model $(W(K),\rho,\tau)$, where   $G=\{x\in W(K)\mid \rho(x)=1\}$.

\begin{example}\label{shyper}{\rm Let $\sign:=\{-1,0,1\}$, ($1+1 = 1,~1-1 = \sign$) be the hyperfield of signs. The associated $W$-model  is $(\Q,\rho,\tau)$ where 
the morphism $\rho:\Q\to \sign$ is given by the sign of rational numbers (\cf\cite{CC5}). The corresponding $\sss$-algebra is then $\cA=H\Q/\Q_+^\times$. One has  $\cA(1_+)=\sign$ and $\cA(2_+)$ is the set of half lines $L$ from the origin in the rational plane $\Q^2$ (including the degenerate case $L=\{0\}$). The maps $\cA(\alpha):\cA(2_+)\to \cA(1_+)$ and $\cA(\beta):\cA(2_+)\to \cA(1_+)$ (\cf~\eqref{alphabeta}) are the projections onto the two axes. The map $\cA(\gamma):\cA(2_+)\to \cA(1_+)$ (\cf~\eqref{addition}) is the projection on the main diagonal. It follows that the hyper-operation \eqref{hypersum} gives back the hyperfield structure on $\sign$. 

One may wonder how to relate the $\sss$-algebra  $\cA=H\Q/\Q_+^\times$ with $H\B$. To this end, one first considers the $\sss$-subalgebra $\cA_+=H\Q_+/\Q_+^\times$ which is defined using the sub-semiring $\Q_+\subset \Q$. The subset $\cA_+(k_+)\subset \cA(k_+)$ corresponds to the collection of half lines $L$ through the origin 
in $\Q^k$ which belong to the first quadrant $\Q_+^k$. One then defines a morphism of $\sss$-algebras $H\rho:H\Q_+\longrightarrow H\B$ by using the morphism of semirings  $\rho:\Q_+\to \B$ and one also notes that $H\rho$ induces a morphism of   $\sss$-algebras $H\rho:\cA_+=H\Q_+/\Q_+^\times\longrightarrow H\B$. One can thus describe the relation between $\cA$ and $H\B$ by the following map
$$
\cA\supset \cA_+\stackrel{H\rho}{\longrightarrow} H\B.
$$
We expect that a similar diagram holds more generally when passing from a semifield of characteristic $1$ to a hyperfield and its
$W$-model.
}\end{example}

\begin{figure}
\begin{center}
\includegraphics[scale=0.33]{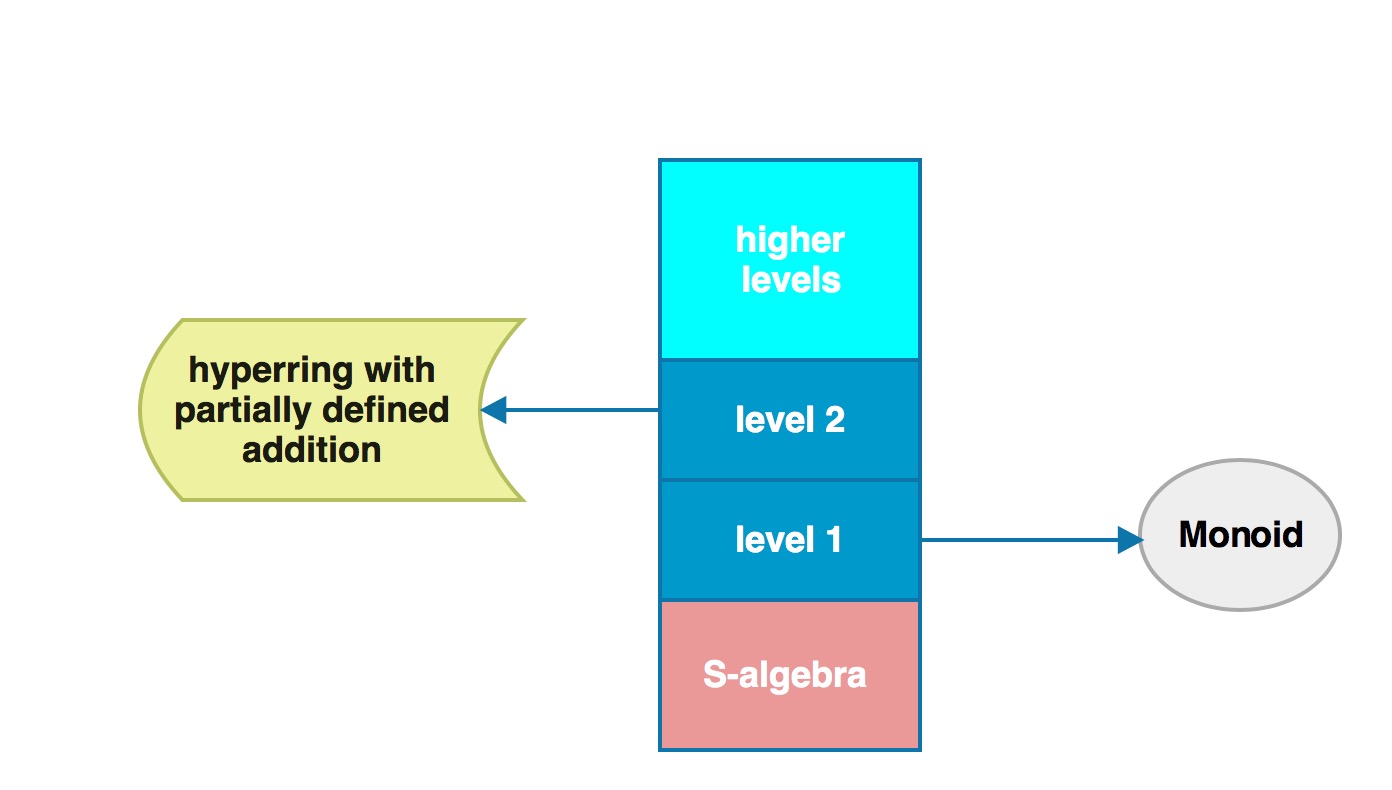}
\end{center}
\caption{Levels of $\sss$-algebras.\label{levels2} }
\end{figure}

\subsection{The levels of an $\sss$-algebra}
In this  section we  briefly discuss the concept of  levels for $\sss$-algebras using which one may derive an explicit approximation  of an $\sss$-algebra $\cA:\gop\longrightarrow \Ses$. The first two levels of such approximation are obtained by restricting the functor $\cA$  to  sets with $\leq 2$ elements to  define the level $1$, and to sets with $\leq 3$ elements to produce the level $2$. In this way, we obtain\vspace{.05in}

{\bf Level $1$}~~The theory of multiplicative mono\"{\i}ds

{\bf Level $2$}~~The theory of hyperrings with partially defined addition.\vspace{.05in}

In \S \ref{sectlevel1} we have seen that at level $1$ the $\sss$-algebra $\cA$ is well approximated by the $\sss$-subalgebra corresponding to the mono\"{\i}d $\cA(1_+)$.
Next we discuss the compatibility of the additive structure given in \eqref{hypersum} with the morphism 
of $\sss$-algebras $\rho:\sss \cA(1_+)\longrightarrow \cA$. For a general mono\" id $M$, implementing \eqref{hypersum} to $\sss M$ with the same notations used there, one gets  
$$
x \oplus y:=\{\sss M(\gamma)z\mid z\in \sss M(2_+), \  \sss M(\alpha)z=x, \ \sss M(\beta)z=y\}\qqq x,y\in \sss M(1_+).
$$
One has: $\sss M(1_+)=M$, while the elements of the set $\sss M(2_+)=2_+\wedge M$  are the base point $0$ and the pairs $(j,m)$ for $j=1,2$, $m\in M$, $m\neq 0$. One also has 
$$
\sss M(\alpha)(1,m)=0, \ \ \sss M(\beta)(1,m)=m, \ \ \sss M(\alpha)(2,m)=m, \ \ \sss M(\beta)(2,m)=0 
$$
and $\sss M(\gamma)(1,m)=\sss M(\gamma)(2,m)=m$. Thus the only rule related to the addition that  is retained at level $1$ simply states that the base point $0$ plays the role of the neutral element:
\begin{equation}
0+x=x+0=x \qqq x\in \sss M(1_+).
\end{equation}
This rule is of course preserved by the morphism 
of $\sss$-algebras $\rho:\sss \cA(1_+)\longrightarrow \cA$.
In \S \ref{fromrtohyp} we have seen that the notion of an $\sss$-algebra is compatible with the operation of quotient
by a subgroup of the multiplicative mono\"{\i}d associated to level $1$. In the following \S \ref{sectnorm} we shall  describe
how to associate a sub-$\sss$-algebra to a sub-multiplicative seminorm on a semiring. By combining these two constructions one provides a good approximation to the description of level $2$ of an $\sss$-algebra. 
This viewpoint shows that, in general, the operation defined by \eqref{hypersum} yields a partially defined hyper-addition so that the approximation of level $2$ is by hyperrings with partially defined addition. The theory of $\sss$-algebras contains however all (non-negative) levels, even though only for the first  level we have an easy explicit description.

\section{$\sss$-algebras and  Arakelov geometry}\label{sectnorm}

This section is motivated by the theory developed in \cite{Durov}  aiming to determine an ``absolute'', algebraic foundation underlying Arakelov geometry: in particular, we refer  to the theory of monads which defines the groundwork of \opcit First, in \S \ref{sectsemi} we explain  how the basic construction of \opcit with the seminorms can be adapted in the framework of $\sss$-algebras. Then, in \S \ref{sectstructsh} we describe  a natural structure sheaf $\cO_\spzb$ of $\sss$-algebras associated to the Arakelov compactification $\spzb$ of $\spz$. As a topological space, $\spzb$ is obtained from the Zariski topological space $\spz$  by adding a new point $\infty$ treated as an additional closed point. This compactification is the same topological space  as the space denoted by $\widehat{\Spec\Z}$ in \opcit (\cf~\S~\ref{durov1}). We show that the Arakelov divisors $D$ on  $\spzb$ give rise to sheaves of $\cO_\spzb$-modules and that the smash product rule for these modules parallels the composition rule of Frobenius correspondences introduced in  \cite{CCas}.

\subsection{$\sss$-algebras and  seminorms}\label{sectsemi}

Let $R$ be a semiring. By a sub-multiplicative seminorm on $R$ we mean a map $R\to\R_+$,  $x\mapsto \vvert x\vvert$, such that $\vvert 0\vvert=0$, $\vvert 1\vvert=1$ and the following inequalities hold
 \begin{equation}\label{subverti1}
\vvert x+y\vvert \leq \vvert x\vvert+\vvert y\vvert, \  \  \vvert xy\vvert \leq\vvert x\vvert\vvert y\vvert \qqq x,y \in R.
\end{equation} 
  The basic construction of \cite{Durov} adapts directly as follows

\begin{prop}\label{sssalg2} $(i)$~Let $R$ be a semiring and $\vvert\cdot \vvert$ a sub-multiplicative seminorm on $R$. Then the following set-up defines an $\sss$-subalgebra $\vvert HR\vvert_1\subset HR$ 
\begin{equation}\label{subvert}
\vvert HR\vvert_1: \Gamma^{{\rm op}}\longrightarrow \Se_*\qquad \vvert HR\vvert_1(X):=\{\phi\in HR(X)\mid \sum_{X\setminus \{*\}} \vvert\phi(x)\vvert\leq 1\}.
\end{equation}
$(ii)$~Let $E$ be an $R$-module and $\vvert\cdot \vvert^E$ a seminorm on $E$ such that $\vvert a\xi\vvert^E\leq \vvert a\vvert \vvert \xi\vvert^E$, $\forall a\in R$, $\forall \xi \in E$. Then for any real $\lambda>0$ the following structure defines a module $\vvert HE\vvert^E_\lambda$ over $\vvert HR\vvert_1$
\begin{equation}\label{subvert1}
\vvert HE\vvert^E_\lambda: \Gamma^{{\rm op}}\longrightarrow \Se_*\qquad \vvert HE\vvert^E_\lambda(X):=\{\phi\in HE(X)\mid \sum_{X\setminus \{*\}} \vvert\phi(x)\vvert^E\leq \lambda\},
\end{equation}
where $HE$ is defined as in \eqref{functsum}.
\end{prop}
\proof $(i)$~The sub-multiplicative seminorm on $R$ fulfills \eqref{subverti1}.
In particular, the triangle inequality shows that, using the formula $HR(f)(\phi)(y)=\sum_{x\mid f(x)=y} \phi(x)$, one derives
$$
\sum_{y\in Y,\,y\neq *}\vvert HR(f)(\phi)(y)\vvert\leq \sum_{y\in Y,\,y\neq *}\left( \sum_{x\mid f(x)=y}\vvert \phi(x)\vvert\right) \leq\sum_{X\setminus \{*\}} \vvert\phi(x)\vvert\leq 1.
$$
Thus, for any morphism $f\in \Hom_{\fin_*}(X,Y)$, the map $HR(f)$ restricts to a map $\vvert HR\vvert_1(X)\to \vvert HR\vvert_1(Y)$. Finally, the stability under product is derived from the sub-multiplicativity of the seminorm as follows
$$
\sum_{(X\wedge Y)\setminus \{*\}} \vvert\phi\psi(z)\vvert=\sum_{X\setminus \{*\}}\sum_{Y\setminus \{*\}}\vvert\phi(x)\psi(y)\vvert\leq\left(\sum_{X\setminus \{*\}}\vvert\phi(x)\vvert\right)\left(\sum_{Y\setminus \{*\}}\vvert\psi(y)\vvert\right)\leq 1.
$$
$(ii)$~The proof of $(i)$ shows that $\vvert HE\vvert^E_\lambda$ is a $\Gamma$-set. It also shows that 
one obtains a map of sets, natural in both objects $X$, $Y$ of $\gop$ 
$$
\vvert HR\vvert_1(X)\wedge \vvert HE\vvert^E_\lambda(Y)\to \vvert HE\vvert^E_\lambda(X \wedge Y), \  \  (a\wedge \xi)(x,y):=a(x)\xi(y), \  \forall x\in X,\ y\in Y.
$$
This gives the required $\vvert HR\vvert_1$-module structure $\vvert HR\vvert_1\wedge \vvert HE\vvert^E_\lambda\longrightarrow \vvert HE\vvert^E_\lambda$ on $\vvert HE\vvert^E_\lambda$.
\endproof

\begin{rem}\label{ss1}{\rm $(i)$~The only multiplicative seminorm on the semifield $\B$ is given by: $\vvert 0\vvert=0$, $\vvert 1\vvert=1$. In this case one obtains $\vvert H\B\vvert_1=\sss$ since the condition \eqref{subvert} restricts the functor $P: \Gamma^{{\rm op}}\longrightarrow \Se_*$ of Lemma \ref{sssalg1} to the range of the unit map $1:\sss\longrightarrow H\B$.\newline
$(ii)$~With the notations of \eqref{subvert1}, one also defines a submodule $\vvert HE\vvert^E_{<\lambda}\subset \vvert HE\vvert^E_\lambda$ by setting
\begin{equation}\label{subvert1strict}
\vvert HE\vvert^E_{<\lambda}: \Gamma^{{\rm op}}\longrightarrow \Se_*\qquad \vvert HE\vvert^E_{<\lambda}(X):=\{\phi\in HE(X)\mid \sum_{X\setminus \{*\}} \vvert\phi(x)\vvert^E< \lambda\}.
\end{equation}
}\end{rem}

\subsection{The structure sheaf  of $\spzb$}\label{sectstructsh}

We let  $\spzb$ be the set $\spz \cup \{\infty\}$ endowed with the topology whose non-empty open sets $U$ are the complements of finite sets $F$ of primes, with $F$ possibly containing $\infty$. In particular all non-empty open sets contain the generic point $(0)\in \spz$.  For any open set $U\subset \spzb$ the subset $U \cap \spz\subset \spz$ is open. The structure sheaf $\cO_\spz$ is canonically a subsheaf of the constant sheaf $\Q$: it associates to the complement of a finite set $F$ of primes the ring $\cO_\spz(F^c)$ of fractions with denominators only involving primes in $F$. We view this ring as a subring of $\Q$.  We shall construct the structure sheaf $\cO_\spzb$ as a subsheaf of the constant sheaf of $\sss$-algebras $H\Q$. First, we use the functor $H$ to transform $\cO_\spz$   into a sheaf of  sub-$\sss$-algebras
$H\cO_\spz\subset H\Q$.  Then, we extend  this sheaf at the archimedean place $\infty\in \spzb$ by 
associating to any open set $U$ containing $\infty$ the $\sss$-subalgebra of $H\Q$ given by
\begin{equation}\label{subsalg}
\hq(U):=\vvert H\cO_\spz(U\cap \spz)\vvert_1
 \end{equation}
where $U\cap \spz=U\setminus\{\infty\}$, $\vvert\cdot\vvert$ is the usual absolute value on $\Q$ which induces a multiplicative seminorm on $\cO_\spz(U\cap \spz)\subset \Q$ so that \eqref{subsalg} makes sense using \eqref{subvert}.\newline
We recall that an Arakelov divisor on $\spzb$ is given by a pair $D=(D_\fii,D_\infty)$ where $D_\fii$ is a divisor on $\spz$ and $D_\infty$ is  a real number $a\in \R$. One uses the notation $D_\infty=a\{\infty\}$. The principal divisor $(q)$ associated to a rational number $q\in \Q^\times$ is $(q):=(q)_\fii -\log\vert q\vert \{\infty\}$. To a divisor $D_\fii$ on $\spz$ corresponds a sheaf $\cO(D_\fii)$ of modules over $\cO_\spz$ and this sheaf is canonically a subsheaf of the constant sheaf $\Q$.
\begin{prop}\label{propspecz} 
$(i)$~The definition \eqref{subsalg} together with the inclusions $\hq(U)\subset H\Q$ and the restriction maps define a sheaf $\hq$ of $\sss$-algebras over $\overline{\Spec\Z}$.

$(ii)$~Let $D$ be an Arakelov divisor on $\overline{\Spec\Z}$. The following defines a sheaf $\cO(D)$ of $\hq$-modules over $\overline{\Spec\Z}$ extending the sheaf $\cO(D_\fii)$ on $\Spec\Z$
\begin{equation}\label{sheafsmod}
\cO(D)(U):=\vvert H\cO(D_\fii)(U\cap \spz)\vvert_{e^a}, \  \ D_\infty=a\{\infty\}\qqq U\ni \infty.
 \end{equation}
The sheaf $\cO(D_\fii)$ on $\Spec\Z$ is viewed as a subsheaf of the constant sheaf $\Q$ which gives a meaning to $\vvert\cdot\vvert$.
\end{prop}
\proof $(i)$~By construction the $\sss$-algebras $\hq(U)$ are naturally $\sss$-subalgebras of $H\Q$. Moreover the restriction map $\hq(U)\to \hq(U')$ for $U'\subset U$ corresponds to the inclusion as $\sss$-subalgebras of $H\Q$. These maps are injective and compatible with the product. To show that one has a sheaf of $\sss$-algebras it is enough to show that for each $k\geq 0$, the (pointed) sets $\hq(U)(k_+)$ form a sheaf of sets. Given a covering $U=\cup U_i$, the elements $\phi:k_+\to \Q$  of $\hq(U_i)(k_+)$ which agree on the intersections $U_i\cap U_j$ correspond to the same element of $H\Q(k_+)$. Moreover this element belongs to $\hq(U)(k_+)$ since this means that  $\phi(j)\in \cO_\spz(U\cap \spz)$ for all $\ell\in k_+$, $j\neq 0$ and, if $\infty \in U$, that $\sum_{\ell\neq *} \vert \phi(\ell)\vert\leq 1$ which is true since $\infty\in U_j$ for some $j$.

$(ii)$~Let $U$ be a non-empty open set containing $\infty$. We apply Proposition \ref{sssalg2} with $R=\cO_\spz(U\cap \spz)$ and the module $E=\cO(D_\fii)(U\cap \spz)$: these are naturally subsets of $\Q$ with the inherited ring and module structures. Thus the norms  induced   by the usual norm on $\Q$ satisfy the hypothesis of Proposition \ref{sssalg2} and one obtains using the same argument as in $(i)$, a sheaf $\cO(D)$ of $\hq$-modules over $\spzb$. Note that at level $1$, one has, for $\infty \in U$,
$$
\cO(D)(U)(1_+)=\{q\in \cO(D_\fii)(U\cap \spz)\mid \vert q\vert \leq e^a\}
=\{q\in\Q^\times\mid D+(q)\geq 0 \ \text{on} \ U\}
$$ 
where the positivity of a divisor is defined pointwise.\endproof

Let $\cA$ be an $\sss$-algebra, $M$ a right $\cA$-module and $N$ a left $\cA$-module.
By definition (\cf~\cite{DGM} Definition 2.1.5.3) the smash product $M\wedge_\cA N$ is the coequalizer 
$$
M\wedge_\cA N:=\varinjlim \{ M\wedge \cA\wedge N \rightrightarrows M\wedge N\}
$$ 
where the maps come from the two actions. We shall use the following two properties of the smash product: the first states that
\begin{equation}\label{propsmash1}
M\wedge_\cA \cA=M, \  \  \cA\wedge_\cA N=N
\end{equation}
and the second asserts that for filtered colimits of $\cA$-modules one has
\begin{equation}\label{propsmash2}
\varinjlim M\wedge_\cA N_k= M\wedge_\cA\varinjlim N_k, \qquad  \varinjlim M_k\wedge_\cA N=\varinjlim M_k\wedge_\cA N.
\end{equation}
These properties can be shown by applying,  for any $\sss$-module $P$, the  isomorphism
$
\Hom_\sss(M\wedge_\cA N,P)= \Hom_\cA(N, \underline{\Hom}(M,P))
$,
where $\underline{\Hom}$ is the internal one.
Likewise for the tensor product of sheaves  modules over sheaves of rings, the smash product $\cO(D)\wedge_\hq \cO(D')$ is defined as the sheaf associated to the presheaf  
$$
U\mapsto \cO(D)(U)\wedge_{\hq(U)}\cO(D')(U).
$$
A similar argument as in Proposition \ref{propspecz} then applies showing that one has a natural morphism of $\hq$-modules 
\begin{equation}\label{multi}
m(D,D'):\cO(D)\wedge_\hq \cO(D')\to \cO(D+D').
\end{equation}\vspace{.05in}

Next, we prove that the smash product rule for the sheaves $\cO(D)$, when taken over the structure sheaf $\hq$ is similar to the law of composition of Frobenius correspondences over the arithmetic site of \cite{CCas}. Let $D_\infty=a\{\infty\}$ and $D'_\infty=a'\{\infty\}$. Then the corresponding conditions at the level $1$ (of the $\sss$-algebras) determine the restrictions 
$
\{q\in\Q^\times\mid \vert q\vert \leq e^a\}
$ and $
\{q\in\Q^\times\mid \vert q\vert \leq e^{a'}\}
$.
 The additive law, asserting that $m(D,D')$ is an isomorphism, 
 means here that the map $(q,q')\mapsto qq'$ is surjective, and fails to be so when both $\lambda=e^a$ and $\lambda'=e^{a'}$ are irrational while their product is rational. The sheaf $\hqr$, which plays the role of the infinitesimal deformation of the identity correspondence in \cite{CCas}, is defined as follows (using \eqref{subvert1strict})
\begin{equation}\label{subsalgr}
\hqr(U):=\vvert H\cO(U\cap \spz)\vvert_{<1}.
 \end{equation}

\begin{prop}\label{proppic} 
$(i)$~The  multiplication rule \begin{equation}\label{multi0}
\cO(D)\wedge_\hq \cO(D')\simeq \cO(D+D')
\end{equation}
holds for the $\hq$-modules $\cO(D)$ and $\cO(D')$ 
except when $\lambda=e^a$ and $\lambda'=e^{a'}$  are irrationals while their product is rational. In that case one has 
\begin{equation}\label{multi0bis}
\cO(D)\wedge_\hq \cO(D')\simeq \hqr \wedge_\hq \cO(D+D')
\simeq  \cO(D+D') \wedge_\hq \hqr.
\end{equation}
$(ii)$~Two Arakelov divisors $D,D'$ are equivalent modulo principal divisors if and only if the sheaves $\cO(D)$ and $\cO(D')$  are isomorphic as abstract sheaves of $\hq$-modules.
\end{prop}
\proof $(i)$~One checks that $m(D,D')$ is an isomorphism, when both $\lambda=e^a$ and $\lambda'=e^{a'}$ are rational: this is verified using \eqref{propsmash1} and the fact that in that case 
the modules are locally trivial. One then obtains \eqref{multi0} in the case $\lambda'=e^{a'}$ is rational using  an increasing sequence $a_n\to a$ such that $e^{a_n}\in\Z[1/p]$ for a fixed prime $p$. On the open set $U=\spzb\setminus \{p\}$ one gets that $\cO(D)(U)=\varinjlim \cO(D_n)(U)$ where $D_n$ is defined by replacing the component at $\infty$ of $D$ by $a_n\{\infty\}$. One then applies   \eqref{propsmash2} to obtain \eqref{multi0}. Finally, when $\lambda'=e^{a'}$ is also irrational one uses  an increasing sequence $a'_n\to a'$ such that $e^{a'_n}\in\Z[1/p]$
together with \eqref{propsmash2} to obtain $(i)$.

$(ii)$~For any $q\in \Q^\times$, one has a canonical isomorphism $\cO(D)\to \cO(D+(q))$ of sheaves of $\hq$-modules over $\spzb$ which is defined by multiplication by $q^{-1}$. More precisely, the module $E=\cO(D_\fii)(U\cap \spz)$ maps to $E'=\cO(D_\fii+(q)_\fii)(U\cap \spz)$ by letting $\xi\mapsto q^{-1}\xi$. One also has  
$$
\sum_{X\setminus \{*\}} \vvert\xi(x)\vvert\leq \lambda\iff 
\sum_{X\setminus \{*\}} \vvert q^{-1}\xi(x)\vvert\leq \lambda \vert q\vert^{-1}.
$$ 
This shows that if $D,D'$ are equivalent modulo principal divisors, then the sheaves $\cO(D)$ and $\cO(D')$ are isomorphic as abstract sheaves of $\hq$-modules.  Conversely,   an isomorphism of  abstract sheaves of $\hq$-modules induces an isomorphism on the set of global sections $\cO(D)(\spzb)(1_+)$. Thus one obtains an invariant of the isomorphism class of abstract sheaf by setting 
$$
\alpha(\cO(D)):=\inf \{b\in \R\mid (\cO(D)\wedge_\hq \cO(b\{\infty\}))(\spzb)(1_+)\neq \{0\}\}
$$
where $D_b:=b\{\infty\}$ is the Arakelov divisor with vanishing finite part. One also has the identification $\cO(D_b)(\spzb)(1_+) = \{n\in\Z\mid \vert n\vert \leq e^b\}$ thus 
$$
\cO(D_b)(\spzb)(1_+)\neq \{0\}\iff b\geq 0.
$$
Using the multiplication rule $(i)$ one gets that  $\alpha(\cO(D_a))=\inf \{b\in \R\mid a+b>0\}=-a$. Since every Arakelov divisor is equivalent, modulo principal divisors, to a $D_a$ one gets $(ii)$.
\endproof

\begin{rem} \label{picrem}{\rm To handle  the non-invertibility of the sheaf $\cO(D)$ for irrational values of $e^a$ in \eqref{sheafsmod} one can use $\hqr$ (\cf~\eqref{subsalgr}) rather than $\hq$ as the structure sheaf.  
The obtained $\sss$-algebras $\cO_\epsilon(U)$ are no longer unital but one still retains the equality $\cO_\epsilon(U)\otimes_{\cO_\epsilon(U)}\cO_\epsilon(U)=\cO_\epsilon(U)$. One defines in a similar manner the sheaves $\cO_<(D)$, using the strict inequality. The multiplication rule is now
\begin{equation}\label{multi1}
\cO_<(D)\wedge_\hqr \cO_<(D')\simeq \cO_<(D+D')
\end{equation}
 This implies that the group of isomorphism classes of invertible sheaves of $\hqr$-modules  contains at least the quotient of the idele class group of $\Q$ by the maximal compact subgroup $\hatz$. Thus this new feature bypasses the defect of the construction of \cite{Durov} mentioned in the open question 4. \S 6 of \cite{durov}.}\end{rem}

\begin{rem} \label{cyclicrem}{\rm The above developments suggest that one should investigate  the (topological) cyclic homology of $\overline{\Spec\Z}$, viewed as a topological space endowed with a sheaf of $\sss$-algebras  and apply the topological Hochschild theory to treat the cohomology of the modules $\cO(D)$.}\end{rem}

\begin{rem}\label{ss12}{\rm Proposition \ref{sssalg2}, combined with the general theory of schemes developed by B. T\"oen and M. Vaqui\'e in \cite{TV} for a symmetric closed monoidal category, provides the tools to perform the same constructions as in \cite{Durov} (reviewed briefly below in \S \ref{durov1}), in the framework of $\sss$-algebras. We shall not pursue this venue, rather we will now explain, using the assembly map of \cite{Lyd}, a conceptual way to pass from the set-up of monads to the framework of $\sss$-algebras.
}\end{rem}

\section{The assembly map}
The main  fact highlighted in this section is that  the assembly map of \cite{Lyd} determines a functorial way to associate an $\sss$-algebra to a monad on $\Ses$. We derive as a consequence  that several main objects introduced in \cite{Durov} can be naturally incorporated in the context of $\sss$-algebras. We start off in \S \ref{durov1} by reviewing some basic structures defined in \opcit   In \S \ref{sectassemb} we  recall the definition of the assembly map which plays  a fundamental role, in  Proposition \ref{assembly}, to define the functor that associates an $\sss$-algebra to a monad.  Finally, Proposition \ref{tensq} of  \S \ref{sectalgmon} shows that the smash product  $\sss$-algebra $H\Z\wedge H\Z$ is (non-trivial and) not isomorphic to $H\Z$. This result is in sharp contrast with the statement of \cite{Durov}  that $\Z \otimes_{\F_1}\Z\simeq \Z$ (\cf~5.1.22 p. 226).
\subsection{Monads and $\widehat{\Spec\Z}$}\label{durov1}
A monad on $\Se$ is defined by an endofunctor $\Sigma:\Se\longrightarrow \Se$ together with two natural transformations $\mu:\Sigma\circ \Sigma\to \Sigma$ and $\epsilon: \id\longrightarrow \Sigma$ which are required to fulfill certain coherence conditions (\!\!\cite{Godem, MacL}).
The required associativity of the product $\mu$ is encoded by the commutativity of the following diagram, for any object $X$ of $\Se$
\begin{equation}\label{1o}
\xymatrix{
\Sigma^3(X) \ar[d]^{\mu_{\Sigma(X)}}  \ar[rr]^{\Sigma(\mu(X))}&& \ar[d]^{\mu_{X}}  \Sigma^2(X) \\
\Sigma^2(X) \ar[rr]_{\mu_X}&&\Sigma(X)
}
\end{equation}
There is a similar compatibility requirement for the unit transformation $\epsilon$. This set-up corresponds precisely to the notion of a mono\"{\i}d in the monoidal category $(\End\,\Se, \circ ,\id)$ of endofunctors under composition, where $\id$ is the identity endofunctor. Next, we  review briefly some of the constructions of \cite{Durov}.
Let $R$ be a semiring, then in \opcit one encodes $R$ as the endofunctor of $R$-linearization $\Sigma_R:\Se \longrightarrow \Se$. It associates to a set $X$ the set  $\Sigma_R(X)$ of finite, formal linear combinations of points of $X$, namely  of finite formal sums $\sum_j \lambda_j\cdot  x_j$, where 
$\lambda_j\in R$ and $x_j\in X$. One uses the addition in $R$ to simplify: $\lambda\cdot x+\lambda'\cdot x=(\lambda+\lambda')\cdot x$. To a map of sets $f:X\to Y$ one associates the transformation 
$$
\Sigma_R(f): \Sigma_R(X)\to \Sigma_R(Y), \qquad  \Sigma_R(f)\left(\sum_j \lambda_j\cdot  x_j\right):=\sum_j \lambda_j\cdot f(x_j).
$$
The multiplication in $R$ is encoded by the elementary fact that a linear combination of finite linear combinations is still a finite linear combination, thus one derives
\begin{equation}\label{mu}
\mu\left(\sum_i \lambda_i\cdot\left(\sum_j \mu_{ij}\cdot  x_{ij}\right) \right)= \sum_{i,j}  \lambda_i \mu_{ij} \cdot x_{ij}.
\end{equation}
This construction defines a natural transformation $\mu:\Sigma_R\circ \Sigma_R\longrightarrow \Sigma_R$ which fulfills the algebraic rules of a monad. 
In fact, the rules describing a monad allow one to impose certain restrictions on the possible linear combinations which one wants to consider. This fits well, in particular, with the structure of the convex sets so that one may for instance define the monad $\Z_\infty$  by considering  as the ring $R$  the field $\R$ of real numbers and by restricting to consider only the linear combinations $\sum \lambda_j x_j$ satisfying the convexity condition $\sum \vert \lambda_j \vert \leq 1$. One can eventually replace $\R$ by $\Q$ and in this case one obtains the monad  $\Z_{(\infty)}$. One may combine the monad $\Z_{(\infty)}$ with the ring $\Z$ as follows. The first step takes an integer $N$ and considers the localized ring $B_N:=\Z[\frac 1 N]$ and the intersection $A_N=B_N\cap \Z_\infty$. This is simply the monad $A_N$ of finite linear combinations  $\sum \lambda_j x_j$ with $\lambda_j\in B_N$ and $\sum \vert \lambda_j \vert \leq 1$.
By computing its prime spectrum one finds that  $\Spec\,A_N=\{\infty\}\cup \Spec \Z\setminus \{p\mid p\vert N\}$. One then reinstalls the missing primes \ie the set $\{p\in \Spec \Z\mid  p\vert N\}$, by gluing $\Spec\,A_N$ with $\Spec \Z$ on the large open set that they have in common. One obtains in this way the set
$$
\widehat{\Spec\Z}^{(N)}:=\Spec\,A_N\cup_{\Spec B_N}\Spec \Z.
$$
Finally, one eliminates the integer $N$ by taking a suitable projective limit under divisibility. The resulting space is the projective limit space
$$
\widehat{\Spec\Z}:=\varprojlim \widehat{\Spec\Z}^{(N)}.
$$
As a topological space it has the same topology (Zariski) as that of $\Spec \Z$ with simply one further closed point (\ie $\infty$) added. At the remaining non-archimedean primes this space fulfills the same topological properties as $\Spec \Z$. One nice feature of this construction is that the local algebraic structure at $\infty$ is modeled by $\Z_{(\infty)}$.

\subsection{The assembly map}\label{sectassemb}
As explained in \S 2.2.1.1 of  \cite{DGM},  a $\Gamma$-set $M:\gop\longrightarrow \Ses$ automatically extends, using 
filtered colimits over the finite subsets $Y\subset X$, to a pointed endofunctor 
\begin{equation}\label{endo1}
\tilde M:\Ses\longrightarrow \Ses, \  \  \tilde M(X):=\varinjlim_{Y\subset X} M(Y).
\end{equation}
 Let $R$ be a (commutative) semiring and $HR$ be the associated $\sss$-algebra as in Lemma \ref{sssalg}. Then the  extension $\widetilde{HR}:\Ses\longrightarrow \Ses$  is, except for the nuance given by the implementation of the base point in the construction, the same as
 the endofunctor $\Sigma_R:\Se\longrightarrow\Se$ of \S \ref{durov1}.

Next, we  review the definition of the assembly map in the simple case of $\Gamma$-sets: we refer to \cite{Lyd, Schw} for more details.   Let  $N:\gop\longrightarrow \Ses$ be a $\Gamma$-set, and $\tilde N: \Ses\longrightarrow \Ses$ its extension to a pointed endofunctor  as in \eqref{endo1}.  One derives a  map natural in the two pointed sets $X,Y$ (\cf~ \S 2.2.1.2 of \cite{DGM})
$$
\iota^N_{X,Y}: X\wedge \tilde N(Y)\to \tilde N(X\wedge Y) 
$$
which is obtained as the family, indexed by the elements $x\in X$, of maps $\tilde N(\delta_x):\tilde N(Y)\to \tilde N(X\wedge Y)$ with 
$$
\delta_x\in \Hom_\Ses(Y,X\wedge Y), \    \   \delta_x(y):=(x,y)\in X\wedge Y \qqq y\in Y.
$$
Similarly, we let  $M:\gop\longrightarrow \Ses$ be a second $\Gamma$-set, and $\tilde M$ its extension to a pointed endofunctor of $\Ses$. One then constructs a  map natural in $X$ and $Y$ as follows
$$
\alpha_{X,Y}:\tilde M( X)\wedge \tilde N(Y)\to \tilde M(\tilde N(X\wedge Y) ),\qquad  \alpha_{X,Y}:=\tilde M(\iota^N_{X,Y})\circ \iota^M_{\tilde N(Y),X}
$$
where  $\iota^M_{\tilde N(Y),X}:\tilde M(X)\wedge \tilde N(Y)\to \tilde M(X\wedge \tilde N(Y))$ and    
$\tilde M(\iota^N_{X,Y}):\tilde M(X\wedge \tilde N(Y))\to \tilde M(\tilde N(X\wedge Y))$.
By restricting to finite pointed sets $X,Y$, \ie to objects of $\gop\subset \fin_*\subset \Ses$, one obtains a natural transformation of bi-functors 
and hence a morphism of $\Gamma$-sets 
\begin{equation}\label{assemblymap}
\alpha_{M,N} :M\wedge N\longrightarrow \tilde M\circ N.
\end{equation}
\begin{defn}\label{assmap} Let  $M,N:\gop\longrightarrow \Ses$ be two  $\Gamma$-sets. The assembly map 
$\alpha_{M,N} :M\wedge N\longrightarrow \tilde M\circ N$ is the morphism of $\Gamma$-sets as in \eqref{assemblymap}.
\end{defn}
\begin{prop}\label{assembly0}
$(i)$~Let $R$ be a semiring, then the assembly map $\alpha_{HR,HR}$
is surjective. 

$(ii)$~For $R=\B$, the assembly map $\alpha_{H\B,H\B}$  associates to a $k$-relation $C=(F,G,v)$  the set of subsets of $k_+$ defined by 
$
\alpha_{H\B,H\B}(C)=\{\{v(x,y)\mid y\in G\}\mid x\in F\}.
$
\end{prop}
\proof $(i)$~The functor $\widetilde{HR}$ coincides with the endofunctor  $\Sigma_R:\Ses \longrightarrow \Ses$. For a pointed set $Z$, we denote the elements of the set $\widetilde{HR}(Z)$  as formal sums 
$\sum r_z\cdot z$, where only finitely many non-zero terms occur in the sum, and the coefficient of the base point $\star$ is irrelevant, \ie the term $r\cdot \star$ drops out. Let $
\xi=\sum r_x\cdot x\in \widetilde{HR}(X)$ and $\eta=\sum s_y\cdot y\in \widetilde{HR}(Y)$, then one has 
$$
\iota^{HR}_{\widetilde{ HR}(Y),X}(\xi \wedge \eta)=\sum r_x\cdot (x \wedge \eta)\in \widetilde{HR}(X\wedge \widetilde{HR}(Y)).
$$
The image of $x \wedge \eta$ by  $\iota^{\widetilde{HR}}_{X,Y}$ is described by $
\iota^{\widetilde{HR}}_{X,Y}(x \wedge \eta)=\sum s_y\cdot (x\wedge y)\in \widetilde{HR}(X\wedge Y)$
and one finally obtains
\begin{equation}\label{assemb1}
\alpha_{HR,HR}(\xi \wedge \eta)=\sum r_x\cdot \left(\sum s_y\cdot (x\wedge y)\right)\in \widetilde{HR}( \widetilde{HR}(X\wedge Y)).
\end{equation}
By applying \eqref{colim} one has 
\begin{equation}\label{colim3}
(HR\wedge HR)(Z)=\varinjlim_{v:X\wedge Y\to Z}( HR(X)\wedge HR(Y)).
\end{equation}
Any element  $\omega\in(HR\wedge HR)(Z)$ is represented by $(X,Y,v,\xi\wedge\eta)$ with $v:X\wedge Y\to Z$ and $\xi, \eta$ as above. Its image is described, using \eqref{assemb1}, by 
\begin{equation}\label{assemb2}
\alpha_{HR,HR}(X,Y,v,\xi\wedge\eta)=\sum r_x\cdot \left(\sum s_y\cdot v(x\wedge y)\right)\in \widetilde{HR}( \widetilde{HR}(Z)).
\end{equation}
For $k_+$ an object of $\gop$, an element of $(\widetilde{HR}\circ HR)(k_+)$  is of the form
$$
\tau=\sum_{i \in I}\lambda_i\cdot\left(\sum_{j\in k_+\setminus \star} \mu_{ij}\cdot  j\right), \   \  \lambda_i\in R, \ \mu_{ij}\in R,
$$
where $I$ is a finite set of indices. Let $X=I_+$ and $\xi=\sum_I \lambda_i\cdot i\in HR(X)$. Let $Y=I_+\wedge k_+$ and $\eta=\sum_{i\in I, \, j\in k_+\setminus \star} \mu_{ij}\cdot  (i,j)\in HR(Y)$. Let $v:X\wedge Y\to k_+$ be given by $v(i,(i',j)):=\star$ if $i\neq i'$ and $v(i,(i,j)):=j$, $\forall i\in I, \, j\in k_+\setminus \star$. With this choice for  $\omega=(X,Y,v,\xi\wedge\eta)$, one then obtains 
$$
\alpha_{HR,HR}(X,Y,v,\xi\wedge\eta)=\sum \lambda_i\cdot \left(\sum \mu_{i'j}\cdot v(i\wedge (i',j))\right)=\tau.
$$
$(ii)$~For $R=\B$,  the element of $(H\B\wedge H\B)(k_+)$ associated to the $k$-relation $C=(F,G,v)$  corresponds to $\omega=(X,Y,v,\xi\wedge\eta)$ where $X=F_+$, $Y=G_+$, $\xi(x)=1$ $\forall x\in F$, $\eta(y)=1$ $\forall y\in G$. Thus \eqref{assemb2} becomes $
\alpha_{H\B,H\B}(C)=\{\{v(x,y)\mid y\in G\}\mid x\in F\}$. 
\endproof

Notice that $(ii)$ of Proposition~\ref{assembly0} implies that the assembly map is {\em not} injective for $R=\B$, since by Corollary \ref{pp1}  one knows that $(H\B\wedge H\B)(k_+)$ is infinite and countable for any $k>0$, while $(\widetilde{H\B}\circ H\B)(k_+)$
remains finite for any choice of $k$.  

\subsection{$\sss$-algebras and monads} \label{sectalgmon}
 The constructions of \cite{Durov} for monads reviewed in  \S \ref{durov1}  can be easily reproduced in the context of $\Gamma$-sets and yield the $\sss$-algebras: $\Z_\infty:=\vvert H\R\vvert_1$ and  $\Z_{(\infty)}:=\vvert H\Q\vvert_1$, 
where the multiplicative seminorms are \resp the usual absolute value on $\R$ and its restriction to $\Q$. This fact  suggests the existence of a functor that associates to a monad on $\Ses$ an $\sss$-algebra. We first briefly comment on  the use of the category of pointed sets in place of the category of sets as in \cite{Durov}. One has a pair of adjoint functors $(L,F)$ where $L:\Se\longrightarrow \Ses$ is the functor $L(X):=X_+$ of implementation of a base point. It is left adjoint to the forgetful functor $F:\Ses\longrightarrow \Se$. The unit of the adjunction is the   natural transformation  $\eta:\id_\Se\longrightarrow F\circ L$ which maps the set $X$ to its copy in the disjoint union $F\circ L(X)=X\,\coprod \,\{\star\}$. The counit  $\epsilon: L\circ F\longrightarrow \id_\Ses$ is the natural transformation which is the identity on $X\subset L\circ F(X)=X\,\coprod\, \{\star\}$ while it maps the extra base point to the base point of $X$.

Let $H$ be a pointed endofunctor of $\Ses$, the composition $H^\#:=F\circ H\circ L$ is then an endofunctor of $\Se$. Moreover given a natural transformation  $\mu:H_1\circ H_2\longrightarrow H_3$ of endofunctors of $\Ses$, one obtains a natural transformation of endofunctors of $\Se$, $ H_1^\#\circ  H_2^\#\longrightarrow  H_3^\#$ as follows
$$
H_1^\#\circ  H_2^\#= F\circ H_1\circ L\circ F\circ H_2\circ L\stackrel{\epsilon}{\longrightarrow} 
 F\circ H_1\circ H_2\circ L\stackrel{\mu}{\longrightarrow} 
F\circ H_3\circ L=H_3^\#.
$$
This construction determines  a natural correspondence between monads on $\Ses$ and monads on $\Se$. Note in particular that the monad on $\Se$ given by $ \id^\#$, where $\id:\Ses\longrightarrow \Ses$ is the identity endofunctor is the monad called $\F_1$ in \cite{Durov}. Moreover, for  any semiring  $R$, with $\widetilde{HR}$ being the extension of $HR$ to an endofunctor of $\Ses$ as in \eqref{endo1}, one checks that the natural monad structure on $\widetilde{HR}$ corresponds to the monad $\Sigma_R$, \ie that ${ \widetilde{HR}^\#}=\Sigma_R$.

Following \opcit, we say that an endofunctor of $\Ses$ is {\em algebraic}  if it preserves filtered colimits: it is then the extension, preserving filtered colimits, of its restriction to $\gop\subset \fin_*\subset \Ses$.
\begin{prop}\label{assembly}
Let $\Sigma$ be a pointed algebraic  monad  on $\Ses$, then the restriction $\cA$ of $\Sigma$ to $\gop\subset \fin_*\subset \Ses$ defines an $\sss$-algebra with the product $m=\mu\circ \alpha_{\Sigma,\cA}$ defined by the composite of the assembly map and the product $\mu:\Sigma\circ \Sigma\to \Sigma$ of the monad $\Sigma$.
\end{prop}
\proof One uses the fact stated in Remark 2.19 of \cite{Lyd} that the assembly map makes the identity functor on $\Gamma$-sets a lax monoidal functor from the monoidal category $(\gam, \circ ,\sss)$ to the monoidal category $(\gam, \wedge,\sss)$. More explicitly, the composition $m=\mu\circ \alpha_{\Sigma,\cA}$ becomes a morphism of $\Gamma$-sets $\cA\wedge \cA\longrightarrow \cA$. The compatibility of the assembly map with the associativity  shows that the commutative diagram \eqref{1o} gives the commutativity of the following diagram of maps of $\Gamma$-sets
\begin{equation}\label{2o}
\xymatrix{
\cA\wedge \cA\wedge \cA \ar[d]^{m\wedge \id}  \ar[rr]^{\id \wedge m}& &\ar[d]^{\mu\circ \alpha_{\Sigma,\cA}}  \cA\wedge \cA \\
\cA\wedge \cA \ar[rr]_{\mu\circ \alpha_{\Sigma,\cA}}& &\cA
}
\end{equation}
The compatibility with the unit  is handled in a similar manner. \endproof 

Using \eqref{assemb1} one checks that, given a semiring $R$ the $\sss$-algebra associated to the monad $\widetilde{HR}\sim \Sigma_R$ is $HR$ (use Proposition \ref{assembly}).

The assembly map makes the identity functor on $\Gamma$-sets a lax monoidal functor from the monoidal category $(\gam, \circ ,\sss)$ to the monoidal category $(\gam, \wedge,\sss)$. Thus given an algebraic monad $\Sigma$ on $\Ses$ one can associate to a (left or right) algebraic module $E$ on $\Sigma$ in the monoidal category $(\gam, \circ ,\sss)$ a (left or right) module over the $\sss$-algebra $\cA$ of Proposition \ref{assembly}. The underlying $\Gamma$-set is unchanged and the action of $\cA$ is obtained using the assembly map as  in Proposition \ref{assembly}. Besides the notion of (left or right) algebraic module $E$ on $\Sigma$  in the sense of the monoidal category $(\gam, \circ ,\sss)$ as discussed in  \S 4.7 of \cite{Durov}, a simpler category $\Sigma-{\rm Mod}$ of 
 modules over a monad is also introduced and used in 3.3.5 of \opcit    A module $E$ on a monad $\Sigma$ on $\Ses$ in this sense is simply a pointed set $E$ together with a map of pointed sets $\alpha:\Sigma(E)\to E$ fulfilling the two conditions: $
\alpha \circ \mu_E=\alpha\circ \Sigma(\alpha), ~\alpha\circ \epsilon_E=\id_E
$.
This notion of a module over a monad defined in \cite{Durov} implies that modules over $\Sigma_R$ correspond exactly to ordinary $R$-modules (\cf~\opcit or \cite{durov} just before Proposition 2, p. 11).  In particular, by taking the monad $\Sigma=\id$ (which is denoted by $\F_\emptyset$ in \cite{Durov}), one finds that any set $M$ is a module in the unique manner provided by the map $\alpha=\id_M:M\to M$. This fact continues to hold if one considers the monad which corresponds to $\F_1$ (\cf ~\cite{durov} \S 4.4 p. 24). Then, one finds that the modules over $\F_1$ are just pointed sets.

 In the framework of $\sss$-algebras, an $R$-module $E$ over a ring $R$ gives rise to an $HR$ module $HE$ but a general module over  $HR$ is not always  of this form (\cf~\cite{DGM} Remark 2.1.5.2, 4). This nuance between $R$-modules and $HR$-modules luckily disappears at the homotopy level. Since this is a crucial fact, we shall comment on it here below: we refer to \cite{Schw1} and \cite{DGM} for the full treatment. Let $k$ be a {\em commutative} $\sss$-algebra, then one defines the  
smash product $M\wedge_k N$ of two $k$-modules as the coequalizer  (\cf~\cite{DGM}, Definition 2.1.5.3)
$$
M\wedge_k N:=\varinjlim \{ M\wedge k\wedge N \rightrightarrows M\wedge N\}
$$ 
where the two maps are given by the actions of $k$ on $M$ and $N$. One obtains in this way the symmetric closed  monoidal category of $k$-modules and a corresponding notion of $k$-algebra. In order to perform homotopy theory, one passes to the associated categories of simplicial objects. The fundamental fact (\cf~\cite{Schw1}, Theorem 4.1) is that the homotopy category of $H\Z$-algebras, obtained using the simplicial version of  $H\Z$-algebras, is equivalent to the homotopy category of ordinary simplicial rings. This result continues to hold if one replaces $\Z$ by any commutative ring $R$ and simplicial rings by simplicial $R$-algebras. This shows that, in the framework of $\sss$-algebras, nothing is lost by using the more relaxed notion of $HR$-module in place of the more restricted notion based on monads.

\subsection{$H\Z\wedge H\Z$ is not isomorphic to $H\Z$}
In \cite{Durov} (\cf~5.1.22 p. 226) one gets the equality $\Z \otimes_{\F_1}\Z\simeq \Z$  which is disappointing for the development of the analogy with the geometric theory over functions fields. The following result shows the different behavior  of the corresponding statement over $\sss$.
\begin{prop}\label{tensq}
The smash product $\sss$-algebra $H\Z\wedge H\Z$ is not isomorphic to $H\Z$.
\end{prop}
\proof Assume that there exists an isomorphism of $\sss$-algebra $\sigma:H\Z\wedge H\Z\stackrel{\sim}{\to} H\Z$. It is then unique since, by Proposition \ref{sssalg1}, the $\sss$-algebra $H\Z$ has no non-trivial automorphism. Furthermore, by Proposition \ref{assembly0} the assembly map $\alpha_{H\Z,H\Z}$
is surjective and this fact implies that the composition $\psi=\alpha_{H\Z,H\Z} \circ  \sigma^{-1}:H\Z\longrightarrow \widetilde{H\Z}\circ H\Z$ is also surjective. Let $\mu:\widetilde{H\Z}\circ H\Z\longrightarrow H\Z$  be the natural transformation as in \eqref{mu}, then the composite $\mu\circ \psi$ is a morphism of $\Gamma$-sets: $\mu\circ \psi\in \Hom_\gop(H\Z,H\Z)$. Since   $\mu:\widetilde{H\Z}\circ H\Z\longrightarrow H\Z$ is surjective, one  then derives that the composite  $\mu\circ \psi\in \Hom_\gop(H\Z,H\Z)$ is also surjective, thus coming from a surjective endomorphism (\ie an automorphism) of the abelian group $\Z$ (here we refer to the proof of 
Proposition \ref{sssalg1}). These facts would imply that $\psi$ is also injective, hence an isomorphism 
$\psi:H\Z\stackrel{\sim}{\to} \widetilde{H\Z}\circ H\Z$. The set $H\Z(k_+)$ is $\Z^k$ with $0$ as the base point. The maps $\alpha,\beta,\gamma:2_+\to 1_+$ of \eqref{addition}  act as follows 
$$
H\Z(\alpha)(n,m)=n, \ \ H\Z(\beta)(n,m)=m, \ \ H\Z(\gamma)(n,m)=n+m.
$$
Thus $\widetilde{H\Z}\circ H\Z$ is naturally isomorphic to the $\Gamma$-set which assigns to $k_+$ the Laurent polynomials in $k$-variables modulo constants $\Z[t^\pm]^{\otimes k}/\Z$. The 
 three maps $\alpha,\beta,\gamma:2_+\to 1_+$ of \eqref{addition}  determine the following maps $\Z[t_1,t_2]\to \Z[t]$
 $$
(\widetilde{H\Z}\circ H\Z)(\alpha)(P)=P(t,1), \  \ (\widetilde{H\Z}\circ H\Z)(\beta)(P)=P(1,t), \  \ (\widetilde{H\Z}\circ H\Z)(\gamma)(P)=P(t,t).
$$
In particular, we see that the map $\rho:(\widetilde{H\Z}\circ H\Z)(2_+)\to (\widetilde{H\Z}\circ H\Z)(1_+)^2$ given by the pair of maps $((\widetilde{H\Z}\circ H\Z)(\alpha),(\widetilde{H\Z}\circ H\Z)(\beta))$ is not injective since one can add to $P(t_1,t_2)$ any multiple of $(t_1-1)(t_2-1)$ without affecting $P(t,1)$ and $P(1,t)$. For the $\Gamma$-set $H\Z$ the map  $\rho':(H\Z)(2_+)\to (H\Z)(1_+)^2$ given by the pair of maps $((H\Z)(\alpha),(H\Z)(\beta))$ is bijective, thus the $\Gamma$-set $(\widetilde{H\Z}\circ H\Z)$ cannot be isomorphic to $H\Z$ and this determines a contradiction.
 \endproof

The elaboration of $H\B\wedge H\B$ seen in \S \ref{sectsmashprod} suggests that the structure of 
$H\Z\wedge H\Z$ is expected to be even far more involved.

\begin{rem}\label{homtop}{\rm When working at the secondary level of topological spectra, it is a well known fact that the smash product 
of the Eilenberg MacLane spectra $\underline{H\Z}\wedge  \underline{H\Z}$ is not isomorphic to $\underline{H\Z}$. In fact, the corresponding homotopy groups are not the same. Indeed, 
for any cofibrant spectrum $X$ one has the equality
$
\pi_*(\underline{H\Z}\wedge X)=H_*(X,\Z)
$,
where $H_*(X,\Z)$ denotes the spectrum homology with integral coefficients. Applying this fact to $X=\underline{H\Z}$,  one interprets $\pi_*(\underline{H\Z}\wedge \underline{H\Z})$ as the integral homology of the Eilenberg-Mac Lane spectrum $\underline{H\Z}$. This homology is known to be finite in positive degrees  but not trivial (\cf~Theorem 3.5 of \cite{Koch}). This argument however, cannot be applied directly to provide an alternative proof of Proposition~\ref{tensq}  since then one would need to compare the spectra 
$
 \underline{H\Z}\wedge  \underline{H\Z}$ and    $   \underline{H\Z\wedge H\Z}
$
and this comparison is usually done by using a cofibrant replacement $H\Z^c$ of the $\Gamma$-set $H\Z$.
}\end{rem}

\end{document}